\documentclass[nonblindrev]{informs3_arxiv}

\usepackage{natbib}
 \bibpunct[, ]{(}{)}{,}{a}{}{,}%
 \def\bibfont{\small}%
\usepackage{amsmath}
\usepackage{amssymb}
\usepackage{multirow}
\usepackage{algorithm}
\usepackage{algorithmic}
\usepackage{setspace}
\usepackage{subfigure}
\usepackage{amsfonts}
\usepackage[table]{xcolor}
\usepackage{color}
\usepackage{graphicx}
\usepackage{epsfig}
\usepackage{enumerate}
\usepackage{epstopdf}
\usepackage{rotating}
\usepackage{bigstrut}
\usepackage{tabularx}
\usepackage{booktabs}
\usepackage{tikz}

\OneAndAHalfSpacedXII 
\newcommand{\sr}{\rule[-0.45cm]{0pt}{0.9cm}}
\newcommand{\st}{\rule[-0.45cm]{0pt}{1.2cm}}
\newcommand{\sn}{\rule[-0.45cm]{0pt}{1cm}}

\usepackage{natbib}
 \bibpunct[, ]{(}{)}{,}{a}{}{,}%
 \def\bibfont{\small}%

\newcommand\hatX{\mathbf{x}^0}

\newcommand\invoptX{\mathbf{x^{*}}}

\newcommand\yinv{\mathbf{y^*}}
\newcommand\cinv{\mathbf{c^*}}

\newcommand{\istar}{ i^* }
\newcommand{\boundary}{\bX^{\textrm{OPT}}}

\newcommand{\bA}{ \mathbf{A} }
\newcommand{\ba}{ \mathbf{a} }
\newcommand{\bbb}{ \mathbf{b} }
\newcommand{\bC}{ \mathbf{C} }
\newcommand{\bc}{ \mathbf{c} }
\newcommand{\bD}{ \mathbf{D} }

\newcommand{\beee}{ \mathbf{e} }

\newcommand{\bp}{ \mathbf{p} }
\newcommand{\bq}{ \mathbf{q} }
\newcommand{\br}{ \mathbf{r} }

\newcommand{\bs}{ \mathbf{s} }
\newcommand{\bt}{ \mathbf{t} }

\newcommand{\bv}{ \mathbf{v} }

\newcommand{\bw}{ \mathbf{w} }
\newcommand{\bX}{ \mathbf{X} }
\newcommand{\bx}{ \mathbf{x} }

\newcommand{\by}{ \mathbf{y} }

\newcommand{\bz}{ \mathbf{z} }

\newcommand{\bzero}{ \mathbf{0} }
\newcommand{\balpha}{ \boldsymbol{\alpha} }
\newcommand{\bxi}{ \boldsymbol{\xi} }

\newcommand{\bepsilon}{ \boldsymbol{\epsilon} }

\newcommand{\blambda}{ \boldsymbol{\lambda} }
\newcommand{\bgamma}{ \boldsymbol{\gamma} }
\newcommand{\bmu}{ \boldsymbol{\mu} }

\DeclareMathOperator*{\sgn}{sgn}

\TheoremsNumberedThrough     
\ECRepeatTheorems
\EquationsNumberedThrough    
\MANUSCRIPTNO{}

\begin{document}

\RUNAUTHOR{T. C. Y. Chan et al.}

\TITLE{Inverse Optimization: Closed-form Solutions, Geometry and Goodness of fit}

\ARTICLEAUTHORS{
\AUTHOR{Timothy C. Y. Chan}
\AFF{Department of Mechanical and Industrial Engineering, University of Toronto, Toronto, Ontario M5S 3G8, Canada, \EMAIL{tcychan@mie.utoronto.ca}} 
\AUTHOR{Taewoo Lee}
\AFF{Department of  Industrial Engineering, University of Houston, Houston, Texas 77204, USA, \EMAIL{tlee6@uh.edu}} 
\AUTHOR{Daria Terekhov}
\AFF{Department of Mechanical and Industrial Engineering, 1455 De Maisonneuve Blvd. W., Montreal, QC H3G 1M8, Canada, \EMAIL{daria.terekhov@concordia.ca}}
}

\ABSTRACT{
In classical inverse linear optimization, one assumes a given solution is a candidate to be optimal. Real data is imperfect and noisy, so there is no guarantee this assumption is satisfied. Inspired by regression, this paper presents a unified framework for cost function estimation in linear optimization comprising a general inverse optimization model and a corresponding goodness-of-fit metric. Although our inverse optimization model is nonconvex, we derive a closed-form solution and present the geometric intuition. Our goodness-of-fit metric, $\rho$, the \emph{coefficient of complementarity}, has similar properties to $R^2$ from regression and is quasiconvex in the input data, leading to an intuitive geometric interpretation. While $\rho$ is computable in polynomial-time, we derive a lower bound that possesses the same properties, is tight for several important model variations, and is even easier to compute. We demonstrate the application of our framework for model estimation and evaluation in production planning and cancer therapy.
}


\maketitle

\section{Introduction}
In inverse optimization, one seeks to impute the cost function or other parameters of an optimization problem in order to render a given solution optimal.  Inverse optimization methods have been developed for linear \citep{Ahuja01}, conic \citep{Iyengar05}, convex \citep{Keshavarz11}, integer \citep{Schaefer09}, multi-objective \citep{Chan14}, variational inequality \citep{Bertsimas14}, and countably infinite linear problems \citep{Ghate15}.  \citet{Heuberger04} surveys inverse network and combinatorial problems.  Inverse optimization has found application in a wide variety of domains, including seismic tomography \citep{Burton92}, demand management \citep{Carr00}, railroad management \citep{Day02}, auctions \citep{Beil03}, production planning \citep{Troutt06}, finance \citep{Bertsimas12}, transportation \citep{Chow12, Chow14}, healthcare \citep{Erkin10, Chan14, Ayer14}, sustainability \citep{Turner13}, and electricity markets \citep{Birge14}.

This paper is motivated by advances to the classical inverse optimization paradigm that arise from data-driven applications. The classical framework implicitly assumes that there exist parameter values that make the given solution optimal (e.g., \citet{Ahuja01}).  However, in many applications it may be impossible to find parameter values that achieve this goal exactly~\citep{Troutt06, Keshavarz11, Chan14, Bertsimas14}. For example, a given datum may be a noisy observation of a pristine solution to an optimization problem, or the model whose parameters are being estimated may be a simpler, lower dimensional model than the one that generated the observation. Applying the classical models in such situations will render the inverse problem infeasible or its solution uninformative (e.g., it may return a trivial solution such as a zero cost vector). Thus, a natural extension of the classical framework is to minimize a measure of error in the fit between model and data in order to make the given solution ``approximately'' optimal. 

A natural way to measure the error in fit between model and data in inverse optimization is via the optimality conditions of the underlying optimization problem. If the optimality conditions are satisfied exactly, then there is no error in the fit. An example of this situation is when the observed solution to a linear optimization problem lies on the boundary of the feasible region -- there exists a cost vector that optimizes that solution exactly. However, if the optimality conditions cannot be satisfied exactly (e.g., the observed solution is an interior point for a linear optimization problem) then there is a positive error in the fit. In the latter case, inverse optimization can identify ``$\epsilon$-optimal'' solutions, i.e., choices of the parameters that enable the observed solution to approximately satisfy the optimality conditions (e.g., \citet{Troutt06,Keshavarz11,Bertsimas14,Chan14}). In the literature, there are several approaches to measuring and minimizing this error in the course of solving an inverse optimization problem. Given its broad applicability, we believe a general measure of model-data fit error in inverse optimization is needed. Note that a few studies have considered context-specific measures of model-data fit error \citep{Troutt06, Chow12}, while others have adapted statistical concepts such as efficiency \citep{Troutt95} and consistency \citep{Aswani15} for inverse optimization.

In this paper, we propose a unified framework for cost function estimation in linear optimization including a general notion of model-data fit error, characterization of the error-minimizing solution to the inverse problem, and a corresponding goodness-of-fit metric.  We draw inspiration from statistics, particularly regression, which has a long history of measuring the fit between a statistical model and real data, intertwined with an optimization model that determines the model parameters.  For example, parameters of a linear regression model are typically derived by minimizing the sum of squared errors between the model and data, which leads to a characterization of the goodness of fit through the coefficient of determination $R^2$.  Just as variants of $R^2$ exist that are tailored to different types of regression approaches (e.g., least squares, least absolute deviation, etc.), we propose variants of our goodness-of-fit metric, $\rho$, that are intimately linked to variations of our general inverse optimization model. We demonstrate the use of our framework to guide model selection in inverse linear optimization. In particular, we use $\rho$ to evaluate and compare models that differ in their constraints on the cost vector. Model selection is a key aspect of developing statistical models and may take the form of constrained regression or variable selection approaches. We demonstrate two analogous concepts for inverse optimization in our numerical results. Since there has been no previous rigorous, general treatment of these topics in inverse optimization, we focus on the case of a single sub-optimal observation as a necessary foundational step toward a framework for the multiple data point case. All proofs are given in the Appendix.

The main contributions of this paper are as follows:

\begin{enumerate}
\item We develop a unified framework for cost function estimation in linear optimization consisting of a general inverse optimization model and a corresponding goodness-of-fit metric. We propose several natural variations of the framework that evaluate goodness of fit in both the decision variables and objective value.

\item We derive closed-form solutions to our general inverse optimization model and to the variants differentiated by the type of model-data error considered.  For each solution, we present the corresponding geometric intuition.

\item We show that our goodness-of-fit metric, $\rho$, the \emph{coefficient of complementarity}, has similar properties to $R^2$ for linear regression: it is maximized by the solution to our inverse optimization model, takes values between 0 and 1, and is non-decreasing in the dimension of the decision vector.  We demonstrate that $\rho$ is quasiconvex in the observed solution, which leads to an intuitive geometric characterization of the goodness of fit within a polyhedron. We also develop an approximation to $\rho$, denoted $\tilde\rho$, that is a lower bound, maintains the same properties, is easier to compute, and coincides with $\rho$ for several special cases.

\item We demonstrate the application of our framework in two domains: production planning and cancer therapy. We use our examples to show how $\rho$ can be used for model selection by evaluating and iteratively improving estimates of a model's parameters via inverse optimization.
\end{enumerate}

\section{Inverse Linear Optimization}\label{sec:inverseOptimization}
In this section, we propose a general inverse optimization model for linear optimization, discuss several natural variants of this model, derive closed-form solutions for each, and illustrate the corresponding geometric intuition.

\subsection{Preliminaries}\label{sec:classicalInverseOptimization}
Let $\bx\in\mathbb{R}^n,\bc\in\mathbb{R}^n,\bA\in\mathbb{R}^{m\times n}$, and $\bbb\in\mathbb{R}^m$.  We define our \emph{forward} optimization problem as %to be the following linear program:
\begin{equation}\label{eq:fop}
\begin{alignedat}{1}
\mathbf{FO}(\bc):\quad \underset{\bx}{\text{minimize}} & \quad \bc'\bx\\
\text{subject to}            & \quad \bA\bx \ge \bbb.
\end{alignedat}
\end{equation}
Let $I = \{1, \dots, m\}$ index the constraint set, $J = \{1, \dots, n\}$ index the set of variables, and $\ba_i$ define the $i$th row of $\bA$. Let $\bX$, assumed full-dimensional and free of redundant constraints, be the set of feasible solutions to $\mathbf{FO}(\bc)$. Let $\bX^{\textrm{OPT}}(\bc)$ be the set of optimal solutions to $\mathbf{FO}(\bc)$, and let $\boundary := \underset{\bc \neq \bzero}\cup \boundary(\bc)$.
We define $\beee_i$ to be a unit vector with the $i$-th component being 1.

Let $\hatX \in \mathbb{R}^n$ denote the \emph{observed solution}. We consider $\bA$, $\bbb$ and $\hatX$ to be exogenously determined; thus, they form the \emph{data} that we use to infer the cost vector. We start with what we refer to as the classical \emph{inverse} optimization problem, which finds a cost vector $\bc$ such that $\hatX$ is optimal for $\mathbf{FO}(\bc)$:
\begin{equation}\label{eqn:IO}
\begin{alignedat}{1}
\mathbf{IO}(\hatX):\quad \underset{\by,\bc}{\text{minimize}} & \quad 0\\
\text{subject to} & \quad \bA'\by = \bc, \\
                  & \quad \bc'\hatX = \bbb'\by, \\
                  & \quad \|\bc\|_1 = 1, \\
                  & \quad \by \ge \bzero.
\end{alignedat}
\end{equation}
We omit $\bA$ and $\bbb$ in our notation $\mathbf{IO}(\hatX)$ for simplicity, but acknowledge that the inverse problem depends on these parameters. The vector $\by$ is the vector of dual variables associated with $\bA\bx \ge \bbb$.  The first two sets of constraints in model~\eqref{eqn:IO} represent dual feasibility and strong duality, respectively.  Besides normalizing the imputed cost vector, the third constraint prevents $\bc = \bzero$ from being a feasible solution to $\mathbf{IO}(\hatX)$.  For problems with nonnegative costs, $\mathbf{IO}(\hatX)$ is a linear program. We omit an objective of the form $\|\bc - \hat\bc\|$, which is present in other approaches (e.g., \citet{Ahuja01,Iyengar05}), that chooses between possibly multiple optimal $\bc$ based on distance to a given $\hat\bc$. In the next section, we will instead propose a model that is flexible enough to accommodate constraints on $\bc$ of the form $\|\bc - \hat{\bc} \| \le \kappa$.

First, we show that for a fixed $\bA$ and $\bbb$, the feasibility of the classical model $\mathbf{IO}(\hatX)$ depends on $\hatX$ in such a way as to severely limit its potential applicability in practice.

\begin{proposition}\label{thm:IO_feasibility}
$\mathbf{IO}(\hatX)$ is feasible if and only if $\hatX \in \bX^{\emph{OPT}}$ or $\hatX \in \{\bx \not\in \bX \;|\; \ba_i'\bx \ge b_i \textrm{ for some } i \in I\}$.
\end{proposition}

The following result is a straightforward corollary (proof omitted) of Proposition~\ref{thm:IO_feasibility}.

\begin{corollary}\label{thm:IO_feasXbounded}
If $\bX$ is bounded, $\mathbf{IO}(\hatX)$ is feasible if and only if $\hatX \in \bX^{\emph{OPT}}$ or $\hatX \not \in \bX$.
\end{corollary}

Proposition~\ref{thm:IO_feasibility} and Corollary~\ref{thm:IO_feasXbounded} highlight the primary undesirable property of $\mathbf{IO}(\hatX)$: a point $\hatX$ that is an interior point of $\bX$ renders $\mathbf{IO}(\hatX)$ infeasible, even if $\hatX$ is very ``close'' to the boundary of $\bX$, in which case there would be at least one natural choice of a cost vector that would make $\hatX$ ``approximately'' optimal.

\subsection{Generalized Inverse Linear Optimization}\label{sec:generalizedInverseOptimization}
We present a simple and natural generalization of $\mathbf{IO}(\hatX)$, where the goal is to simultaneously derive a cost vector $\bc^*$ and identify the ``smallest'' perturbation of $\hatX$ so the perturbed point is in the set $\bX^{\textrm{OPT}}(\bc^*)$. By considering the distance between $\hatX$ and $\bX^{\textrm{OPT}}(\bc^*)$ as a measure of error, we define a natural measure of goodness of fit between the data and the fitted model.
While error in inverse optimization has been explored previously (e.g., \citet{Troutt06,Keshavarz11,Chan14,Bertsimas14}), this paper is the first to leverage such a notion to develop a general goodness-of-fit metric for inverse optimization that is integrated with the underlying inverse optimization model.

Our approach determines a meaningful cost vector for any $\hatX \in \bX$, which is the only restriction we place on $\hatX$ from here on out. Throughout the paper, we view the inverse problem through the lens of fitting, i.e., identifying a cost vector that best fits a (noisy) observation, represented by $\hatX$.

Consider the following model, which generalizes $\mathbf{IO}(\hatX)$:
\begin{equation}\label{eq:GIO_general}
\begin{alignedat}{1}
\mathbf{GIO}(\hatX):\quad \underset{\by,\bc,\bepsilon}{\text{minimize}} & \quad {\|\bepsilon\|_L} \\
\text{subject to} & \quad \bA'\by = \bc,\\
		 & \quad \bc'(\hatX - \bepsilon) = \bbb'\by,\\
                  & \quad \|\bc\|_1 = 1, \\
                  & \quad \by \ge \bzero.
\end{alignedat}
\end{equation}

The objective function $\|\cdot\|_L$ is an arbitrary norm (we use ``$L$'' for ``loss''). The error vector $\bepsilon$ in the strong duality constraint allows $\mathbf{GIO}(\hatX)$ to be feasible when $\hatX$ is an interior point. If $\hatX \in \bX^{\textrm{OPT}}$, then the optimal $\bepsilon$ equals $\bzero$: there is no error between the data and the model. Otherwise, there will be a nonzero $\bepsilon$ representing positive error in the fit.
Technically, a primal feasibility constraint $\bA(\hatX - \bepsilon) \ge \bbb$ should be included in formulation~\eqref{eq:GIO_general}. However, it is automatically satisfied by an optimal solution to $\mathbf{GIO}(\hatX)$, which we will show later. We assume $\bc$ and $\bepsilon$ are unrestricted for now and revisit this issue in Section~\ref{sec:structuralConstraints}. 

Unlike $\mathbf{IO}(\hatX)$, $\mathbf{GIO}(\hatX)$ has an optimal solution given any $\hatX \in \bX$.
\begin{proposition}\label{thm:GIO_feasibility}
Given any $\hatX \in \bX$,
\begin{enumerate}
\item $\mathbf{GIO}(\hatX)$ has an optimal solution.
\item A solution $(\by, \bc, \bepsilon)$ is feasible to $\mathbf{GIO}(\hatX)$ if and only if $(\by, \bc)$ is a feasible solution to $\mathbf{IO}(\hatX - \bepsilon)$.
\end{enumerate}
\end{proposition}

The second part of Proposition~\ref{thm:GIO_feasibility} elucidates the geometry associated with an optimal solution to $\mathbf{GIO}(\hatX)$. Given $\hatX$, $\mathbf{GIO}(\hatX)$ identifies a direction of perturbation, $\bepsilon^*$, of minimal distance (with respect to $\|\cdot\|_L$) to bring $\hatX$ into the set $\bX^{\textrm{OPT}}$. An optimal $\bc^*$ to $\mathbf{GIO}(\hatX)$ is a cost vector such that $\hatX - \bepsilon^* \in \bX^{\textrm{OPT}}(\bc^*)$.  That is, solving $\mathbf{GIO}(\hatX)$ is equivalent to solving $\mathbf{IO}(\bx^*)$ where $\bx^* := \hatX - \bepsilon^*$ is a point in $\bX^{\textrm{OPT}}$ closest to $\hatX$ as measured by $\|\cdot\|_L$.  Thus, the distance from $\hatX$ to $\bX^{\textrm{OPT}}$, $\|\bepsilon^*\|_L$, is a measure of model-data fit error.

\subsection{Characterizing an optimal solution to $\mathbf{GIO}(\hatX)$}\label{sec:GIOwithFeasibleSolution}
Although $\mathbf{GIO}(\hatX)$ is not linear due to the term $\bc'\bepsilon$, we can characterize its optimal solutions in closed form.  As described above, solving $\mathbf{GIO}(\hatX)$ is equivalent to finding a projection of $\hatX$ on to $\bX^{\textrm{OPT}}$ using $\|\cdot\|_L$.  Equivalently, $\hatX$ is being projected on to the closure of the complement of the convex set that contains $\hatX$. The main result of this subsection, the closed-form characterization of optimal solutions to $\mathbf{GIO}(\hatX)$, relies on the concept of a dual norm, which we define here for completeness.  If $\| \cdot \|$ is a norm, then we write the associated dual norm as $\|\bz\|^D := \sup \{\bz'\bx : \|\bx\| \le 1\}$.

\begin{theorem}\label{thm:GGIO}
Given $\hatX \in \bX$, an optimal solution to $\mathbf{GIO}(\hatX)$ is
\begin{equation}\label{eq:closedFormSolution}
(\by^*, \bc^*, \bepsilon^*) = \left(\frac{\beee_{\istar}}{\|\ba_{\istar}\|_1}, \frac{\ba_{\istar}}{\|\ba_{\istar}\|_1}, \frac{\ba_{\istar}'\hatX - b_{\istar}}{\|\ba_{\istar}\|_{L}^{D}} \bv(\ba_{\istar})\right),
\end{equation}
where $\istar \in \underset{i \in I} \argmin \{ (\ba_i'\hatX - b_i)/\|\ba_i\|_L^D\}$ and $\bv(\ba_{\istar}) \in \arg \max_{\|\bv\|_L=1} \ba_{\istar}'\bv$. Furthermore, the optimal objective value of $\mathbf{GIO}(\hatX)$ is $\|\bepsilon^*\|_L = (\ba_{\istar}'\hatX - b_{\istar})/\|\ba_{\istar}\|_L^D$.
\end{theorem}

\begin{remark} Theorem 1 implies that the search for an optimal $\bc$ is reduced to comparing a simple ratio across a finite set of choices and that $\bc^* \in \left\{ \ba_{1}/\|\ba_{1}\|_1, \ba_{2}/\|\ba_{2}\|_1,  \dots, \ba_{m}/\|\ba_{m}\|_1 \right\}$.
\end{remark}

The main takeaway from Theorem~\ref{thm:GGIO} is that a cost vector that makes $\hatX$ minimally suboptimal (according to the $\|\cdot\|_L$ metric) is easily computable: identify the constraint $i$ with the smallest value of $(\ba_i'\hatX - b_i)/\|\ba_i\|_L^D$ and an optimal $\bc^*$ will be a normalized multiple of the vector $\ba_i$ defining that constraint. The fact that the search for an optimal cost vector can be restricted to a finite set of alternatives implies a natural choice for a baseline comparison when measuring goodness of fit; we revisit this issue in Section~\ref{sec:needForGoodnessOfFit} where we formally define our goodness of fit metric, $\rho$. Note that there may exist multiple $i$'s that achieve the $\arg\min$ in Theorem~\ref{thm:GGIO}. In particular, if the point $\hatX$ is projected exactly to a vertex of $\bX$, then a $\bc$ vector that is a suitably normalized conic combination of the corresponding multiple $\ba_i$ vectors will also be optimal.

We now return to an assertion made earlier, namely that the constraint $\bA(\hatX-\bepsilon)\ge \bbb$ is automatically satisfied by an optimal solution $\bepsilon^*$ to $\mathbf{GIO}(\hatX)$. The proof of Theorem~\ref{thm:GGIO} establishes that this assertion is indeed true. While the projection of $\hatX$ to an arbitrary constraint may not be to a feasible point, the projection to the closest constraint measured by $||\cdot||_L$ will be feasible.

\subsection{Variants of $\mathbf{GIO}(\hatX)$}\label{sec:specializations}
In this subsection, we present three natural variants of the general model $\mathbf{GIO}(\hatX)$. These variants arise from different choices of the norm in the objective function of $\mathbf{GIO}(\hatX)$ and possibly restricting the structure of $\bepsilon$, which together induce different projections of $\hatX$ on to $\boundary$. In all three cases, the characterizations of $\bc^*$ and $\by^*$ do not change; only the specific constraint on to which $\hatX$ is projected changes. For completeness, we define the sign function as $\sgn(s) = s/|s|$ for a nonzero $s \in \mathbb{R}$ and 0 otherwise. When applied to a vector $\bs$, the sign function returns a vector of the component-wise signs denoted $\sgn(\bs)$.

\subsubsection{$p$-norm ($\mathbf{GIO}_{p}(\hatX)$)}
First, we consider minimizing error using the $p$-norm, $p \ge 1$, which provides a natural measure of error in the space of decision variables.  Recall that the $p$-norm is defined as $\|\bx\|_p = (\sum_{j=1}^{n} |x_j|^p)^{1/p}$ for $p \in [1, \infty)$ and $\|\bx\|_{\infty} = \underset{1 \le j \le n}\max\{|x_j|\}$ for $p = \infty$.  The resulting formulation sets the objective of formulation~\eqref{eq:GIO_general} to $\|\bepsilon\|_p$:
\begin{equation}\label{eq:GIO_l2}
\begin{alignedat}{1}
\mathbf{GIO}_p(\hatX): \underset{\by,\bc,\bepsilon}{\text{minimize}} & \quad {\|\bepsilon\|}_p \\
\text{subject to} & \quad \bA'\by = \bc,\\
				  & \quad \bc'(\hatX-\bepsilon) = \bbb'\by,\\
                  & \quad \|\bc\|_1 = 1, \\
                  & \quad \by \ge \bzero.
\end{alignedat}
\end{equation}

Solving $\mathbf{GIO}_p(\hatX)$ returns a cost vector $\bc^*$ that finds $\bx^* := \hatX-\bepsilon^* \in \arg\min_{\bx\in\boundary} \|\bx-\hatX\|_p$.  The following corollary specializes Theorem~\ref{thm:GGIO} to the $p$-norm case, with particular attention paid to $p = 1, 2, \infty$ (proof omitted). Recall that the $p$- and $q$-norms are dual if $1/p + 1/q = 1$.  In particular, the 1-norm and $\infty$-norm are dual, and the 2-norm is self-dual. 

\begin{corollary}\label{thm:GIOp_reformulation}
Let $\istar \in \underset{i \in I} \argmin \{ (\ba_i'\hatX - b_i)/\|\ba_i\|_q\}$ where $q$ satisfies $1/p + 1/q = 1$ and $j^* \in \argmax_l|a_{\istar, l}|$.  Then the optimal objective value of $\mathbf{GIO}_p(\hatX)$ is $\|\bepsilon^*\|_p = (\ba_{\istar}'\hatX - b_{\istar})/\|\ba_{\istar}\|_q$. For $p = 1, 2, \infty$, an optimal $\bepsilon^*$ is
\begin{align}
\bepsilon^* & = \left\{
  \begin{array}{l l}
    \sgn(a_{\istar,j^*}) \frac{\ba_{\istar}'\hatX - b_{\istar}}{\|\ba_{\istar}\|_{\infty}}\beee_{j^*}, & \quad \textrm{ if } p = 1, \\
    \frac{\ba_{\istar}}{\|\ba_{\istar}\|_2}\frac{\ba_{\istar}'\hatX - b_{\istar}}{\|\ba_{\istar}\|_2}, & \quad \textrm{ if } p = 2, \\
    \sgn(\ba_{\istar})\frac{\ba_{\istar}'\hatX - b_{\istar}}{\|\ba_{\istar}\|_1}, & \quad \textrm{ if } p = \infty.
  \end{array} \right.
\end{align}
\end{corollary}

Geometrically, $\mathbf{GIO}_p(\hatX)$ finds the largest $p$-norm ball centered at $\hatX$ that lies within the feasible region.  The point of intersection between this ball and the boundary of the feasible region is $\hatX-\bepsilon^*$, which lies on the constraint defined by $\ba_{\istar}$.
\subsubsection{Absolute duality gap ($\mathbf{GIO}_a(\hatX)$)}\label{sec:GIOa}
Next, we consider error measured in terms of the objective function value. In particular, we find $\bc^*$ that minimizes the difference in objective values between $\hatX$ and an optimal solution $\bx^*$ to $\mathbf{FO}(\bc^*)$. This approach can be interpreted as minimizing the absolute duality gap, $\epsilon_a$ \citep{Troutt06,Chan14}:
\begin{equation}\label{eq:ailop_a}
\begin{alignedat}{1}
\mathbf{GIO}_a(\hatX): \underset{\by,\bc,\epsilon_a}{\text{minimize}} & \quad \epsilon_a\\
\text{subject to} & \quad \bA'\by = \bc,\\
                  & \quad \bc'\hatX = \bbb'\by+\epsilon_a,\\
                  & \quad \|\bc\|_1 = 1, \\
                  & \quad \by \ge \bzero.
\end{alignedat}
\end{equation}
An attractive feature of this formulation is that it becomes a linear program when the cost vector is non-negative. Being able to solve the inverse problem as a LP directly may be useful when there are additional application-specific constraints that need to be included (see Sections~\ref{sec:modelSelection} and~\ref{sec:experimentalResults}).

Next, we show that $\mathbf{GIO}_a(\hatX)$ can be derived through a simple modification of $\mathbf{GIO}(\hatX)$.

\begin{proposition}\label{thm:equivalenceGIOa}
$\mathbf{GIO}_a(\hatX)$ is equivalent to $\mathbf{GIO}(\hatX)$ with $\|\cdot\|_L = \|\cdot\|_\infty$ and the constraint $\bepsilon = \epsilon_a \sgn(\bc)$.
\end{proposition}

Even though $\mathbf{GIO}_a(\hatX)$ is not a direct specialization of $\mathbf{GIO}(\hatX)$, it turns out that a result analogous to Theorem~\ref{thm:GGIO} holds.

\begin{proposition}\label{thm:GIOa_reformulation}
An optimal solution to $\mathbf{GIO}_a(\hatX)$ is $(\by^*, \bc^*, \epsilon_a^*) = (\beee_{i^*}/\|\ba_{i^*}\|_1, \ba_{i^*}/\|\ba_{i^*}\|_1, (\ba_{\istar}'\hatX - b_{\istar})/\|\ba_{\istar}\|_1)$, where $\istar \in \underset{i \in I} {\arg\min} \{ (\ba_i'\hatX - b_i)/\|\ba_i\|_1\}$.
\end{proposition}

This result says that when error is measured in terms of the absolute duality gap, a cost vector that minimizes this error will coincide with one of the $m$ constraints $\ba_i$. Geometrically, $\mathbf{GIO}_a(\hatX)$ evaluates level sets parallel to the hyperplanes $\{\bx  |  \ba_i'\bx = b_i\}$ for all $i$ to determine a cost vector $\bc^*$ and a point $\bx^* \in \mathbf{X}^{\textrm{OPT}}(\bc^*)$ that minimizes the absolute duality gap, i.e., ${\bc^*}'\hatX - {\bc^*}'\invoptX$. 

\subsubsection{Relative duality gap ($\mathbf{GIO}_r(\hatX)$)}\label{sec:relative}
A second measure of error in the objective function value is the relative duality gap \citep{Troutt95,Chan14}. Assuming $\bbb'\by \neq 0$, we define the relative duality gap as ${(\bc'\hatX - \bbb'\by)}/{|\bbb'\by|} = |{\bc'\hatX}/\bbb'\by - 1|$. If we define $\epsilon_r := \bc'\hatX/\bbb'\by$, then the inverse model minimizing the relative duality gap is
\begin{equation}\label{eq:GIO_r_formulation}
\begin{alignedat}{1}
\mathbf{GIO}_r(\hatX): \underset{\by,\bc,\epsilon_r}{\text{minimize}} & \quad |\epsilon_r - 1|\\
\text{subject to} & \quad \bA'\by = \bc,\\
                  & \quad \bc'\hatX = \epsilon_r\bbb'\by,\\
                  & \quad \|\bc\|_1 = 1, \\
                  & \quad \by \ge \bzero.
\end{alignedat}
\end{equation}

Like $\mathbf{GIO}_a(\hatX)$, $\mathbf{GIO}_r(\hatX)$ can be derived from a simple modification to $\mathbf{GIO}(\hatX)$.

\begin{proposition}\label{thm:equivalenceGIOr}
$\mathbf{GIO}_r(\hatX)$ is equivalent to $\mathbf{GIO}(\hatX)$ with $\|\cdot\|_L$ equal to the weighted infinity norm $\| \cdot \|_{\infty,K} = \| \cdot \|_{\infty}|K|$, $K = 1/\bbb'\by$, and  the constraint $\bepsilon = \bbb'\by(\epsilon_r - 1) \sgn(\bc)$.
\end{proposition}

Similar to the case of $\mathbf{GIO}_a(\hatX)$, a result analogous to Theorem~\ref{thm:GGIO} holds for $\mathbf{GIO}_r(\hatX)$.

\begin{proposition}\label{thm:GIOr_reformulation}
An optimal solution to $\mathbf{GIO}_r(\hatX)$ is $(\by^*,\bc^*,\epsilon_r^*) = (\beee_{i^*}/\|\ba_{i^*}\|_1, \ba_{i^*}/\|\ba_{i^*}\|_1, \ba_{i^*}' \bx^0 /b_{i^*})$, where $i^* \in \underset{i\in I}{\arg\min} \{(\ba_{i}' \bx^0 - b_i)/|b_i | \}.$
\end{proposition}

This result implies that Remark 1 is again relevant here: there exists an optimal $\bc$ that is a normalized version of one of the $\ba_i$ vectors. The only difference is in the specific ratios that need to be compared to determine the optimal $\bc$. Geometrically, $\mathbf{GIO}_r(\hatX)$ evaluates level sets parallel to the hyperplanes $\{\bx  |  \ba_i'\bx = b_i\}$ $\forall i$ to determine an optimal cost vector $\bc^*$ and a point $\bx^* \in \mathbf{X}^{\textrm{OPT}}(\bc^*)$ that minimizes the relative duality gap, i.e., $|{\bc^*}'\hatX/{\bc^*}'\invoptX - 1|$. With knowledge of the closed-form solution structure, the requirement $\bbb'\by \ne 0$ can be satisfied if at least one of the $b_i$ is not zero.

Finally, we discuss how changing the normalization constraint allows $\mathbf{GIO}_r(\hatX)$ to be solved by solving at most two linear programs. This result may be practically useful when a computational solution to $\mathbf{GIO}_r(\hatX)$ is desired, such as when there are extra constraints like those used for model selection.

\begin{proposition}\label{GIO_r_LP_alternative}
An optimal solution to $\mathbf{GIO}_r(\hatX)$ is $(\by^*,\bc^*,\epsilon_r^*) = (\hat{\by}/\|\hat{\bc}\|_1,\hat{\bc}/\|\hat{\bc}\|_1, \hat{\epsilon}_r)$ where $(\hat{\by},\hat{\bc}, \hat{\epsilon}_r)$ is an optimal solution to formulation~\eqref{eq:GIO_r_alternative}.
\begin{equation}\label{eq:GIO_r_alternative}
\begin{alignedat}{1}
\underset{\by,\bc, \epsilon_r}{\emph{minimize}} & \quad |\epsilon_r - 1|\\
\emph{subject to} & \quad \bA'\by = \bc,\\
                  & \quad \bc'\hatX = \epsilon_r\bbb'\by,\\
                  & \quad |\bbb'\by| = 1, \\
                  & \quad \by \ge \bzero.
\end{alignedat}
\end{equation}
\end{proposition}

Formulation~\eqref{eq:GIO_r_alternative} can be solved by solving two linear programs, one with $\bbb'\by = 1$ and the other with $\bbb'\by = -1$ (and of course linearizing the absolute value in the objective), and then choosing the solution with lowest cost. Furthermore, if $\bbb'\by > 0$, which is the case in many applied problems with positive right hand sides, it is necessary to solve only a single problem \citep{Chan14} and then appropriately rescale its solution:
\begin{equation}\label{eq:GIO_r_single}
\begin{alignedat}{1}
\underset{\by,\bc}{\textrm{minimize}} & \quad \bc'\hatX\\
\textrm{subject to} & \quad \bA'\by = \bc,\\
                  & \quad \bbb'\by = 1, \\
                  & \quad \by \ge \bzero.
\end{alignedat}
\end{equation}

Note that if $\bc \ge \bzero$, formulation~\eqref{eq:GIO_r_single} not only minimizes the relative duality gap, it does so while performing a weighted $\ell_1$ regularization of the cost vector $\bc$ with the weights corresponding to the components of $\bx^0$. This observation suggests that solutions to $\mathbf{GIO}_r(\hatX)$ may offer sparsity, which we revisit in Section~\ref{sec:cancertherapy}. 

\subsection{Summary}\label{sec:summary}
The main results of this section are: 1) determining the optimal $\bc^*$ to $\mathbf{GIO}(\hatX)$ can be done by restricting consideration to a finite set of $m$ cost vectors and comparing $m$ corresponding ratios (Theorem~\ref{thm:GGIO}); 2) variants of $\mathbf{GIO}(\hatX)$ can specialize measurement of error to the space of decision variables or objective function values, while maintaining the simple closed-form optimal solution structure. It is important to note that for the two duality gap models all points on the constraint $\ba_{i^*}'\bx = b_{i^*}$ induce the same error, regardless of whether the point is primal feasible, since error is measured in terms of objective value with respect to $\bc^*$. This property simplifies the computation of the goodness of fit with respect to the objective value, as shown in Section~\ref{sec:goodnessOfFit}.

Table \ref{tab:formulations} summarizes the variants of $\mathbf{GIO}(\hatX)$ presented in the previous subsections. The first column indicates the space in which error is measured.  The next three columns identify the type of model along with the corresponding choices of the objective function norm and assumed structure of $\bepsilon$ used to connect the variant and $\mathbf{GIO}(\hatX)$. The last two columns specify the structure of the closed-form solutions, which can be found by substituting the entries in those columns into equation~\eqref{eq:closedFormSolution}. 

\begin{table}[ht!]
\begin{center}
\begin{tabular}{cccccccc} \hline
& \multicolumn{4}{c}{Model Variant} && \multicolumn{2}{c}{Solution Structure} \\ \cmidrule(l){2-5}\cmidrule(l){7-8}
\sn & \multicolumn{2}{c}{Type} & \multicolumn{1}{c} {$\|\bepsilon\|_L$} &
$\bepsilon$ &&  \multicolumn{1}{c}{$\|\ba_{\istar}\|_L^D$} & $\bv(\ba_{\istar})$ \\ \cmidrule(l){2-3}\cmidrule(l){4-4}\cmidrule(l){5-5}\cmidrule(l){7-7}\cmidrule(l){8-8}
\multicolumn{1}{c}{ \multirow{3}{*}[-3mm]{\begin{sideways}\small Decision space\end{sideways}} } &
\multicolumn{1}{c}{ \multirow{3}{*}[-5mm]{$\mathbf{GIO}_p(\hatX)$} } &
\sr 1 &
$\| \bepsilon \|_1$ &
$\bepsilon$ && $\| \ba_{\istar} \|_{\infty}$
& $\sgn(a_{\istar,j^*})\beee_{j^*}$, $j^*$ from Cor.~\ref{thm:GIOp_reformulation}\\
 \multicolumn{1}{c}{} &
 &
 \sr 2 &
 $\| \bepsilon \|_2$
 & $\bepsilon$
 & & $\| \ba_{\istar} \|_2$
 & $\ba_{\istar}/\|\ba_{\istar}\|_2$  \\ 
 \multicolumn{1}{c}{} &
&  \sr $\infty$
& $\| \bepsilon \|_{\infty}$
& $\bepsilon$
& & $\| \ba_{\istar} \|_1$
& $\sgn(\ba_{\istar})$ \\  
\multicolumn{1}{c}{  \multirow{2}{*}[3mm]{\begin{sideways} \small Objective space\end{sideways}}} &
 \multicolumn{2}{c}{\st $\mathbf{GIO}_a(\hatX)$} & $\| \bepsilon \|_{\infty}$ &
 $\epsilon_a\sgn(\bc)$ & &
 $\| \ba_{\istar} \|_1$ &
 $\sgn(\ba_{\istar})$ \\ 
\multicolumn{1}{c} {} &
\multicolumn{2}{c}{\st $\mathbf{GIO}_r(\hatX)$} &$\| \bepsilon \|_{\infty, 1/|\bbb'\by|}$ & $\bbb'\by(\epsilon_r-1)\sgn(\bc)$ & & $|b_{\istar}|$  & $|b_{\istar}|\sgn(\ba_{\istar})/\|\ba_{\istar}\|_1$ \\
\hline
\end{tabular}
\caption{Summary of the $\mathbf{GIO}(\hatX)$ variants.}
\label{tab:formulations}
\end{center}
\end{table}

The following numerical example illustrates the different geometric interpretations of the $\mathbf{GIO}(\hatX)$ variants.

\begin{figure}
\centering
\includegraphics[scale=0.7]{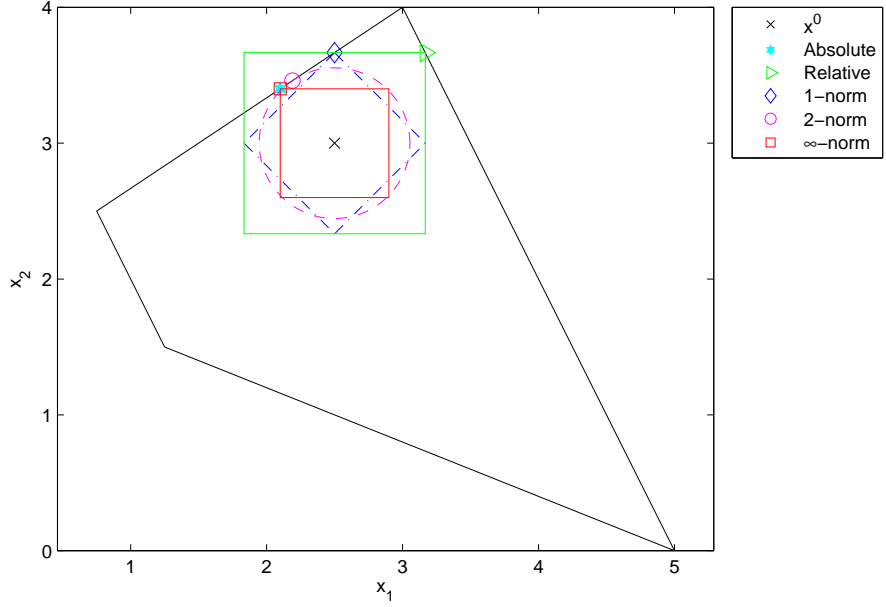}
\caption{$\mathbf{GIO}(\hatX)$ projections for Example 1.}
\label{fig:epsilon_a_1}
\end{figure}

\begin{example}\label{example}
\emph{Consider the following forward problem:
\begin{equation}\label{eq:example}
\begin{alignedat}{1}
\underset{\bx}{\emph{minimize}} & \quad \quad c_1x_1 + c_2x_2 \\
\emph{subject to} & \quad \quad 2x_1 + 5x_2 \ge 10, \\
& \quad \quad 2x_1 - 3x_2 \ge -6, \\
& \quad \quad \; 2x_1 + x_2 \ge 4, \\
& \quad -2x_1 - x_2 \ge -10.
\end{alignedat}
\end{equation}
The feasible region of this problem is shown in Figure \ref{fig:epsilon_a_1}.  Let $\hatX = (2.5, 3)$. The projected solutions $\hatX - \bepsilon^*$ are $[2.5, 3.\bar{6}]$ for $\mathbf{GIO}_1(\hatX)$, $[2.19, 3.46]$ for $\mathbf{GIO}_2(\hatX)$, $[2.1, 3.4]$ for $\mathbf{GIO}_{\infty}(\hatX)$ and $\mathbf{GIO}_a(\hatX)$, and $[3.1\bar{6}, 3.\bar{6}]$ for $\mathbf{GIO}_r(\hatX)$. For $\mathbf{GIO}_p(\hatX)$, the point $\hatX - \bepsilon^*$ is where the largest $p$-norm ball meets the boundary of the feasible region. For $\mathbf{GIO}_a(\hatX)$, the projection to the boundary matches $\mathbf{GIO}_{\infty}(\hatX)$. For $\mathbf{GIO}_r(\hatX)$, the projected point is again at the intersection of a $\infty$-norm ball and the boundary. The specific ball is the one with the smallest weighted norm among the four that just meet each of the four constraints, weighted by the corresponding $b_i$. 
}
\end{example}

\section{Goodness of Fit}\label{sec:goodnessOfFit}

We first briefly review the concept of goodness of fit in linear regression and then develop our goodness-of-fit metric for inverse optimization. 

\subsection{Preliminaries}\label{sec:regressionGoodnessOfFit}

A multiple linear regression model assumes a linear relationship between a dependent variable $Y$ and independent variables $X_1, \ldots, X_n$:
\begin{equation}
Y = \beta_0 + \beta_1 X_1 + \cdots + \beta_n X_n + \epsilon,\label{eqn:regression}
\end{equation}
where $\epsilon$ is a zero-mean, normally distributed error term.  Given data $(\bx^1, y_1), \ldots, (\bx^Q, y_Q)$, the most common method for determining the regression parameters $\beta_0, \ldots, \beta_n$ is to minimize the sum of squared errors or residuals, also known as \emph{linear least squares}:
\begin{equation}
\underset{\beta_0, \ldots, \beta_n}{\textrm{minimize}} \quad \sum_{q=1}^Q \left(y_q - \beta_0 - \sum_{j=1}^n \beta_j x^q_j\right)^2. \label{eqn:leastSquares}
\end{equation}

Correspondingly, the most popular and simple-to-calculate metric to quantify the goodness of fit of a linear regression model is the \emph{coefficient of determination}, $R^2$.  The coefficient of determination can be thought of as the proportion of the total variation in the observed $y_q$ values that is accounted for by the independent variables (i.e., the model).  More precisely,
\begin{align}
R^2 = 1 - \frac{\textrm{variation not explained by model}}{\textrm{total variation}} = 1 - \frac{\sum_{q=1}^{Q} (y_q - \hat{y}_q)^2}{\sum_{q=1}^{Q} (y_q - \bar{y})^2} \label{eqn:R2},
\end{align}
where $\hat{y}_q$ is the value of $Y$ estimated from the model using $\bx^q$ and $\bar{y}$ is the average of the observed $y_q$ values. The coefficient of determination possesses several attractive properties: it is maximized when $\boldsymbol{\beta}$ is determined by minimizing the sum of squared residuals; it takes values in $[0, 1]$ and is therefore easy to interpret; and it is non-decreasing in the number of independent variables in the regression model.

Although $R^2$ as defined in \eqref{eqn:R2} is the most well-known goodness-of-fit metric in linear regression, alternative $R^2$-type metrics exist \citep{Sprecher94}. For example, there is a variant of $R^2$ for least absolute deviation regression that has absolute values instead of squares in both the numerator and denominator of equation~\eqref{eqn:R2}, and replaces the average, $\bar{y}$, with the median of the observed $y_q$ values. This variant of $R^2$ is maximized when $\boldsymbol{\beta}$ is determined by minimizing the sum of absolute values of the residuals, and also takes values in $[0, 1]$ \citep{Pynnonen94}. 

\subsection{Goodness of Fit in Generalized Inverse Optimization}\label{sec:needForGoodnessOfFit}

Inspired by the coefficient of determination in regression, we develop a goodness-of-fit metric for inverse optimization that we call the \emph{coefficient of complementarity}. The name is derived from the satisfaction of the complementary slackness/optimality conditions in $\mathbf{GIO}(\hatX)$. Goodness of fit in inverse optimization has been examined previously in very specific contexts \citep{Troutt06,Chow12}. The metric we present is general, is appropriately linked with the underlying inverse optimization model, and possesses several attractive properties that mirror properties of $R^2$.  We first present our proposed (``exact'') metric, and discuss its interpretation and computation. Then, we discuss an alternative metric that approximates the exact one -- and in some cases coincides with the exact one -- but is easier to compute.

We define the coefficient of complementarity as: 
\begin{equation}\label{eq:singleObsRhoGeneralGIO}
\rho = 1 - \frac{\|\bepsilon^*\|_L}{\frac{1}{m}\sum_{i=1}^m \|\bepsilon^i\|_L},
\end{equation}
where $\|\bepsilon^*\|_{L}$ quantifies the error in data-model fit determined by solving $\mathbf{GIO}(\hatX)$ (i.e., the shortest distance from $\hatX$ to $\boundary$ as measured by $\|\cdot\|_L$), and $\|\bepsilon^i\|_{L}$ measures the error in the fit with respect to constraint $i$, $i=1, \ldots, m$. Note that the dependence of $\rho$ on the given point $\hatX$ is implicit in~\eqref{eq:singleObsRhoGeneralGIO}. The denominator of $\rho$ captures the average error associated with $m$ possible inverse solutions, one for each constraint of the primal feasible region (recall Theorem~\ref{thm:GGIO}). Implicit in the above definition is the assumption that $\sum_{i=1}^m \|\bepsilon^i\|_L \neq 0$. Indeed, if $\sum_{i=1}^m \|\bepsilon^i\|_L = 0$, then $\|\bepsilon^i\|_L = 0$ for all $i$, in which case we define $\rho := 1$.

Intrinsic in the definition of $\rho$ (and $R^2$) is the fact that it is not an absolute measure of performance. Rather, it evaluates the model-data fit of the estimated model (i.e., one that identifies a cost vector that minimizes the error $\|\bepsilon\|_L$) relative to the average variation in the feasible region with respect to $\hatX$, i.e., all the given data in the problem, $(\hatX, \bA, \bbb)$. Alternatively, the denominator of $\rho$ can be viewed as the expected value of a geometric ``null hypothesis'' that $\bc$ is chosen uniformly at random from the set $\{\ba_i \; | \; i \in I \}$. Thus, $\rho$ also captures the value of the inversely optimized solution against such a randomly chosen $\bc$. This null hypothesis is the first attempt to define a general baseline solution that is analogous to $\bar{y}$ in regression.

One of the main advantages of a relative measure like $\rho$ is that it automatically puts in context whether the ``raw'' error, $\|\bepsilon\|_L$, is big or small.
As a unitless measure, $\rho$ is invariant to scaling in $\bA$ and $\bbb$, which is similar to how $R^2$ is invariant to scaling in the input data of a least-squares regression model.
While domain experts working on a specific application may know the acceptable level of raw error, our metric can be used in any application setting by a nonexpert.

Recalling the $p$-norm, absolute duality gap, and relative duality gap variants of $\mathbf{GIO}(\hatX)$ described in Section~\ref{sec:specializations}, $\rho$ also has corresponding variants that measure the appropriate error for each model. While $\|\bepsilon^*\|_L$ can be obtained by Theorem~\ref{thm:GGIO}, calculation of the denominator differs between the variations of $\rho$ based on different geometric interpretations of model-data fit.

\begin{figure}[ht]\centering
\includegraphics[height = 8cm, width=10cm]{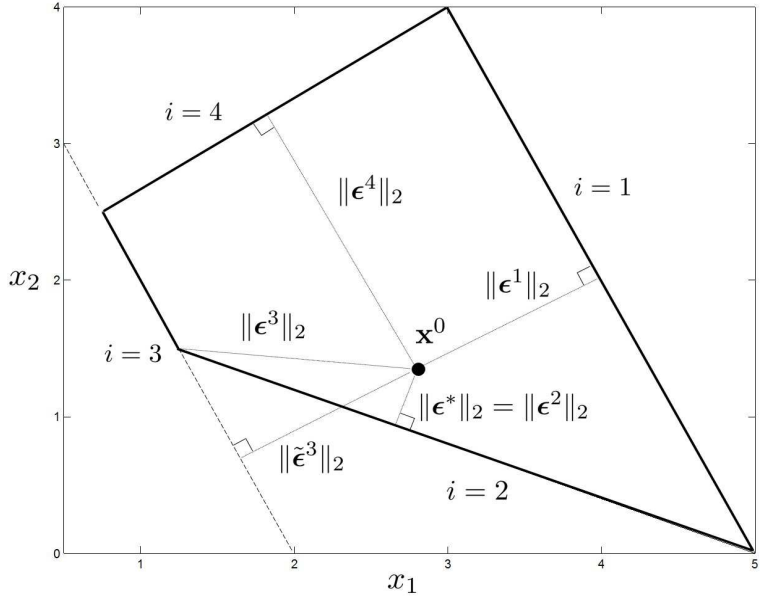}
\caption{Illustration of $\|\bepsilon^*\|_2$, $\|\bepsilon^i\|_2$, and $\|\tilde\bepsilon^i\|_2$ for Example 1.}
\label{fig:rho_i_illustration}
\end{figure}

\subsubsection{$\rho$ for $\mathbf{GIO}_p(\hatX)$.}

We define $\rho_p$, the variant of $\rho$ that corresponds to $\mathbf{GIO}_p(\hatX)$, as
\begin{equation}\label{eq:singleObsRhoGeneralGIO_p}
\rho_p
= 1 - \frac{\|\bepsilon^*\|_p}{\frac{1}{m}\sum_{i=1}^m \|\bepsilon^i\|_p}
= 1 - \frac{\|\bx^* - \hatX\|_p}{\frac{1}{m}\sum_{i=1}^m \|\bx^i-\hatX\|_p},
\end{equation}
where $\bx^* := \hatX - \bepsilon^*$ and $\bx^i := \hatX - \bepsilon^i$. The first expression is a direct substitution of the $p$-norm into equation~\eqref{eq:singleObsRhoGeneralGIO},
while the second emphasizes the fact that error is being measured in the space of decision variables.

To find $\bepsilon^i$, i.e., the error in the fit with respect to constraint $i= 1, \ldots, m$, the following convex optimization problem can be solved:
\begin{equation}\label{eq:eps_p^i}
\begin{alignedat}{1}
\underset{\bepsilon}{\text{minimize}} & \quad \|\bepsilon\|_p\\
\textrm{subject to } & \quad \bA (\hatX - \bepsilon) \ge \bbb, \\
& \quad \ba_i' (\hatX - \bepsilon) = b_i.
\end{alignedat}
\end{equation}

It is important to recognize that $\|\bepsilon^{i}\|_L$ does not, in general, represent the minimum distance from $\bx^0$ to constraint $i$, due to the primal feasibility constraint. Figure~\ref{fig:rho_i_illustration} illustrates $\bepsilon^*$ and $\bepsilon^i$ in the polyhedron from Example~\ref{example} when $\|\cdot\|_L$ is the 2-norm; the minimum distance projection to constraint $i$ is represented by $\tilde\bepsilon^i$, which is an optimal solution to~\eqref{eq:eps_p^i} without the primal feasibility constraint $\bA(\hatX-\bepsilon)\ge\bbb$.

\subsubsection{$\rho$ for $\mathbf{GIO}_a(\hatX)$ and $\mathbf{GIO}_r(\hatX)$.} In the absolute and relative duality gap cases, we know through Propositions~\ref{thm:equivalenceGIOa} and~\ref{thm:equivalenceGIOr} that $\|\bepsilon^*\|_L$ equals $\epsilon^*_a$ and $\epsilon^*_r$, respectively. Similarly, through Propositions~\ref{thm:equivalenceGIOa}--\ref{thm:GIOr_reformulation}, we can define $\|\bepsilon^i\|_L$, the fit error with respect to constraint $i$, as $\epsilon^i_a$ and $\epsilon^i_r$, respectively. Thus, we define $\rho_a$ and $\rho_r$, the variants of $\rho$ that correspond to $\mathbf{GIO}_a(\hatX)$ and $\mathbf{GIO}_r(\hatX)$, respectively, as
\begin{equation}\label{eq:singleObsRhoGeneralGIO_a}
\rho_a = 1 - \frac{\epsilon^*_a}{\frac{1}{m}\sum_{i=1}^m \epsilon^i_a}
 = 1 - \frac{{\bc^*}'\hatX-{\bc^*}'\bx^*}{\frac{1}{m}\sum_{i=1}^m ({\bc^i}'\hatX-{\bc^i}'\bx^i)}
 = 1 -
\frac{(\ba_{i^*}'\hatX-b_{i^*})/\|\ba_{i^*}\|_1}{\frac{1}{m}\sum_{i=1}^m ({\ba_i}'\hatX-b_i)/\|\ba_i\|_1},
\end{equation}
and
\begin{equation}\label{eq:singleObsRhoGeneralGIO_r}
\rho_r = 1 - \frac{|\epsilon^*_r - 1|}{\frac{1}{m}\sum_{i=1}^m |\epsilon^i_r - 1|} = 1 - \frac{|{\bc^*}'\hatX/{\bc^*}'\bx^* - 1|}{\frac{1}{m}\sum_{i=1}^m |{\bc^i}'\hatX / {\bc^i}'\bx^i - 1|} = 1 - \frac{|{\ba_{i^*}}'\hatX/b_{i^*} - 1|}{\frac{1}{m}\sum_{i=1}^m |{\ba_i}'\hatX / b_i - 1|},
\end{equation}
where $\bx^i$ is an optimal solution to $\mathbf{FOP}(\bc^i)$,  $\bx^* = \bx^{i^*}$, $\bc^i = \ba_i / \|\ba_i\|_1$, and $\bc^* = \bc^{i^*}$. The second expression in each of equations~\eqref{eq:singleObsRhoGeneralGIO_a} and~\eqref{eq:singleObsRhoGeneralGIO_r} shows explicitly that error is measured in terms of objective function value, while the third expression is due to Propositions \ref{thm:GIOa_reformulation} and \ref{thm:GIOr_reformulation}. Note that finding $\epsilon^i_a$ and $\epsilon^i_r$ does not require ensuring primal feasibility as in the $p$-norm case (see Section~\ref{sec:summary}). As a result, $\epsilon^i_a$ and $\epsilon^i_r$ can be determined by simply setting $\epsilon^i_a$  and $\epsilon^i_r$ as follows:
\begin{equation}\label{eq:eps_i_a_def}
\epsilon^i_a = \frac{(\ba_i'\hatX - b_i)}{\|\ba_i\|_1},
\end{equation}
\begin{equation}\label{eq:eps_i_r_def}
\epsilon^i_r =\frac{ \ba_i'\hatX}{b_i},
\end{equation}
where the latter assumes $b_i \neq 0$.

\subsection{Properties of $\rho$}\label{sec:theory}
In this sub-section, we present key properties of $\rho$, as defined in~\eqref{eq:singleObsRhoGeneralGIO}, which hold for all $\rho$ variants presented above. Some additional notation is needed before proceeding. We define $\rho^{(k)} = 1 - {\|\bepsilon^{(k)*}\|_L}/({\sum_{i=1}^m \|\bepsilon^i\|_L}/m)$, where $\bepsilon^{(k)*}$ is an optimal solution for $\bepsilon$ for $\mathbf{GIO}^{(k)}(\hatX)$, the generalized inverse optimization model~\eqref{eq:GIO_general} where the first $k \le n$ (without loss of generality) components of $\bc$ are free and the remaining $n-k$ components are fixed to 0. The denominator of $\rho^{(k)}$ matches the denominator of $\rho$, since the denominator represents the average error induced by the problem data $(\bx^0, \bA, \bbb)$, which is independent of the inverse optimization process. We will also write $\rho(\hatX)$ in place of $\rho$ where it is useful to show the explicit dependence of $\rho$ on the given point $\hatX$.

\begin{theorem}\label{thm:rho_properties}
Let $\hatX \in \bX$.
\begin{enumerate}
\item (Optimality) $\rho$ is maximized by an optimal solution of $\mathbf{GIO}(\hatX)$,
\item (Boundedness) $\rho \in [0,1]$,
\item (Monotonicity) $\rho^{(k)} \le \rho^{(k+1)}$ for $k = 1, \ldots, n-1$,
\item (Quasiconvexity) $\rho(\hatX)$ is a quasiconvex function of $\hatX \in \bX$.
\end{enumerate}
\end{theorem}

The first three properties directly mirror the properties possessed by $R^2$. The first property implies that although $\rho$ can be used to evaluate the goodness of fit of a solution generated by any inverse optimization model, it is maximized by an optimal solution to $\mathbf{GIO}(\hatX)$; the regression analogue is that while $R^2$ as defined in equation~\eqref{eqn:R2} can measure the goodness of fit of a regression model with coefficients generated from any method, it is maximized when the coefficients are determined by minimizing the sum of squared residuals (i.e., solving problem~\eqref{eqn:leastSquares}).

The second property shows that $\rho$, just like $R^2$, is a unitless quantity that takes values in $[0, 1]$, with 0 indicating complete lack of fit and 1 indicating perfect fit. Perfect fit corresponds to the existence of a cost vector that makes $\hatX$ optimal for the forward problem.  Indeed, when $\hatX \in \bX^{\textrm{OPT}}$, $\bepsilon^* = \bzero$, so $\rho = 1$. On the other hand, if $\hatX \in \bX \setminus \bX^{\textrm{OPT}}$, then $\bepsilon^* \ne \bzero$, so $\rho < 1$. The case of $\rho = 0$ indicates a complete lack of fit, which occurs when $\|\bepsilon^i\|_L = \|\bepsilon^*\|_L$ for all $i = 1, \ldots, m$, as in the case of $\bX$ being a unit hypercube centered at $\hatX$, for example.

The monotonicity property is the analogue of $R^2$ being nondecreasing in the number of independent variables. This property reflects the fact that having a more flexible model (with more degrees of freedom) enables better fit with the data. However, just like in regression, it may be possible to overfit an inverse optimization model. Thus, the monotonicity property requires the user to be cognizant of artificially increasing $\rho$ by including more variables/dimensions without actually gaining more insight or producing a model that is generalizable beyond the given data.

\begin{figure}\centering
\includegraphics[width=9cm]{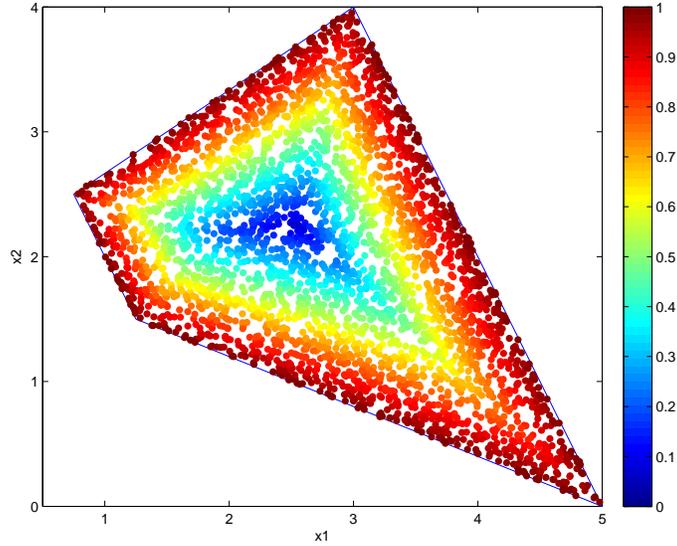}
\caption{$\rho_2$ values throughout the polyhedron from Example~\ref{example}.}
\label{fig:rho_2_illustration}
\end{figure}

The quasiconvexity property implies that $\rho$ varies in a structured way inside the polyhedron $\bX$.
For example,
Figure~\ref{fig:rho_2_illustration} illustrates the $\rho_2$ values for points $\hatX$ scattered throughout the polyhedron from Example~\ref{example}. Since quasiconvexity of a function $f(\bx)$ is equivalent to convexity of its sub-level sets $\{\bx | f(\bx) \le \alpha\}$ for $\alpha \in \mathbb{R}$, geometrically we see that $\rho$ is non-decreasing in a ``radial'' sense as $\hatX$ moves from the interior towards the boundary.
This observation is captured formally in Proposition~\ref{thm:distance_quasiconvexity}.
\begin{proposition}\label{thm:distance_quasiconvexity}
Let $\bx^{\min} \in \{\bx \in \bX \; | \; \rho(\bx) \le \rho(\by) \; \forall \; \by \in \bX\}$ and $\bx^* \in \normalfont{\boundary}$. Let $\bx_1 = \lambda_1 \bx^{\min} + (1-\lambda_1) \bx^*$ and $\bx_2 = \lambda_2 \bx^{\min} + (1-\lambda_2) \bx^*$ where $\lambda_1 \in (0, 1)$, $\lambda_2 \in (0, 1)$, and $\lambda_1 \neq \lambda_2$. If $\|\bx_1 - \bx^*\|_L \le \|\bx_2 - \bx^*\|_L$ then $\rho(\bx_1) \ge \rho(\bx_2)$.
\end{proposition}
Without quasiconvexity, points closer to the boundary may, paradoxically, exhibit worse fit than those that are closer to the ``center'' of the polyhedron. 

Overall, the above-described properties enable $\rho$ to be a general tool to aid model selection in inverse optimization, with its primary attractiveness being its interpretability and applicability by a nonexpert.

\subsection{An alternative $\rho$ metric}
In this subsection, we present an alternative goodness-of-fit metric, $\tilde\rho$, that has essentially the same interpretation and retains the same attractive properties of $\rho$. This metric coincides with $\rho$ in the absolute and relative duality gap cases, and adapts the simple approach of those calculations to the $p$-norm case. Thus, it is simpler to compute than $\rho$ in general.

Analogously to $R^2$, whose computation requires only the solution of the least squares problem~\eqref{eqn:leastSquares} and the given data, the effort required to compute $\rho$ would ideally not extend beyond what is needed to identify $i^*$ and compute $\|\bepsilon^*\|_L$, i.e., to solve the inverse optimization problem.  However, this is not the case for $\rho_p$, which requires the additional effort of solving problem~\eqref{eq:eps_p^i} $m$ times, once for each constraint $i$. Granted, the effort is modest given that~\eqref{eq:eps_p^i} is convex, but it is additional effort nonetheless. While identifying $i^*$ and determining $\|\bepsilon^*\|_L$ simply requires the evaluation of $(\ba_i'\hatX - b_i)/\|\ba_i\|_L^D$ for each $i$, $\|\bepsilon^i\|_L$ is, in general, not equal to $(\ba_i'\hatX - b_i)/\|\ba_i\|_L^D$, since the projection of $\hatX$ on to some constraints may result in a point that is not primal feasible (cf. Figure~\ref{fig:rho_i_illustration}). This observation suggests that a natural approximation for $\rho$, which we denote as $\tilde\rho$, can be derived by simply ignoring the primal feasibility restriction when computing $\bepsilon^i$, which is essentially what the $\rho_a$ and $\rho_r$ computations are doing. 

Formally, we define $\tilde\rho$ as
\begin{equation}\label{eq:approxrho}
\tilde\rho = 1 - \frac{\|\bepsilon^*\|_L}{\frac{1}{m}\sum_{i=1}^m \|\tilde\bepsilon^i\|_L},
\end{equation}
where $\|\tilde\bepsilon^i\|_L = (\ba_i'\hatX - b_i)/\|\ba_i\|_L^D$ is the norm of the projection vector of $\hatX$ on to constraint $i$, ignoring primal feasibility of the projected point (i.e., $\hatX - \tilde\bepsilon^i$ is not required to be primal feasible). For the $p$-norm case, $\tilde\bepsilon^i$ is an optimal solution to~\eqref{eq:eps_p^i} without the constraint $\bA(\hatX-\bepsilon)\ge\bbb$. The computation of $\tilde\rho$ requires only the evaluation of the ratios $(\ba_i'\hatX - b_i)/\|\ba_i\|_L^D$ for each $i$, which is already carried out when solving $\mathbf{GIO}(\hatX)$ (cf. Theorem~\ref{thm:GGIO}). Thus, there is no additional effort required to compute $\tilde\rho$ beyond solving the inverse problem. Note that the only difference between $\tilde\rho$ and $\rho$ is in the denominator. We define $\tilde\rho_p$ analogously to $\rho_p$, except with $\tilde\bx^i := \hatX - \tilde\bepsilon^i$ in place of $\bx^i$ in equations~\eqref{eq:singleObsRhoGeneralGIO_p}. As in the case of $\rho$, we assume that $\sum_{i=1}^m \|\tilde{\bepsilon}^i\|_L \neq 0$; otherwise, we define $\tilde{\rho} := 1$.

First, we show that $\tilde\rho$ retains the same four properties of $\rho$ presented in Section~\ref{sec:theory}. The term $\tilde\rho^{(k)}$ is defined analogously to $\rho^{(k)}$.
\begin{proposition}\label{thm:approxrho_properties}
Let $\hatX \in \bX$.
\begin{enumerate}
\item (Optimality) $\tilde\rho$ is maximized by an optimal solution of $\mathbf{GIO}(\hatX)$,
\item (Boundedness) $\tilde\rho \in [0,1]$,
\item (Monotonicity) $\tilde\rho^{(k)} \le \tilde\rho^{(k+1)}$ for $k = 1, \ldots, n-1$,
\item (Quasiconvexity) $\tilde\rho(\hatX)$ is a quasiconvex function of $\hatX \in \bX$.
\end{enumerate}
\end{proposition}

Next, we examine the relationship between $\tilde\rho$ and $\rho$.  Recall that $\|\tilde\bepsilon^i\|_L$ measures the distance from $\hatX$ to constraint $i$ whereas $\|\bepsilon^i\|_L$ measures the distance from $\hatX$ to the closest feasible point on constraint $i$. So, in terms of the inverse optimization process, using $\tilde\bepsilon^i$ underestimates the lack of data-model fit with respect to constraint $i$ ($\|\tilde\bepsilon^i\|_L \le \|\bepsilon^i\|_L$ for all $i$) and thus $\tilde\rho$ is a lower bound for $\rho$. However, $\tilde\rho$ and $\rho$ coincide when error is measured in terms of objective function value, as in $\mathbf{GIO}_a(\hatX)$ and $\mathbf{GIO}_r(\hatX)$; $\tilde\bepsilon^i$ and $\bepsilon^i$ projecting $\hatX$ to different points on constraint $i$ becomes irrelevant since the cost of those two points will be equal with respect to $\ba_i$. These results are summarized in the following proposition (proof omitted).

\begin{proposition}\label{thm:rhoapprox}
\hfill
\begin{enumerate}
\item $\tilde\rho\le\rho$, 
\item $\tilde\rho = \rho$ for the absolute and relative duality gap variants.
\end{enumerate}
\end{proposition}

\begin{figure}[!ht]\centering
\subfigure[Distribution of $\rho_2$.]
{\label{fig:rho_2_exact_comparison}\includegraphics[width=75mm]{rho_2_plot_0_avg.eps}}
\subfigure[Distribution of $\tilde\rho_2$.]
{\label{fig:rho_2_approx_comparison}\includegraphics[width=75mm]{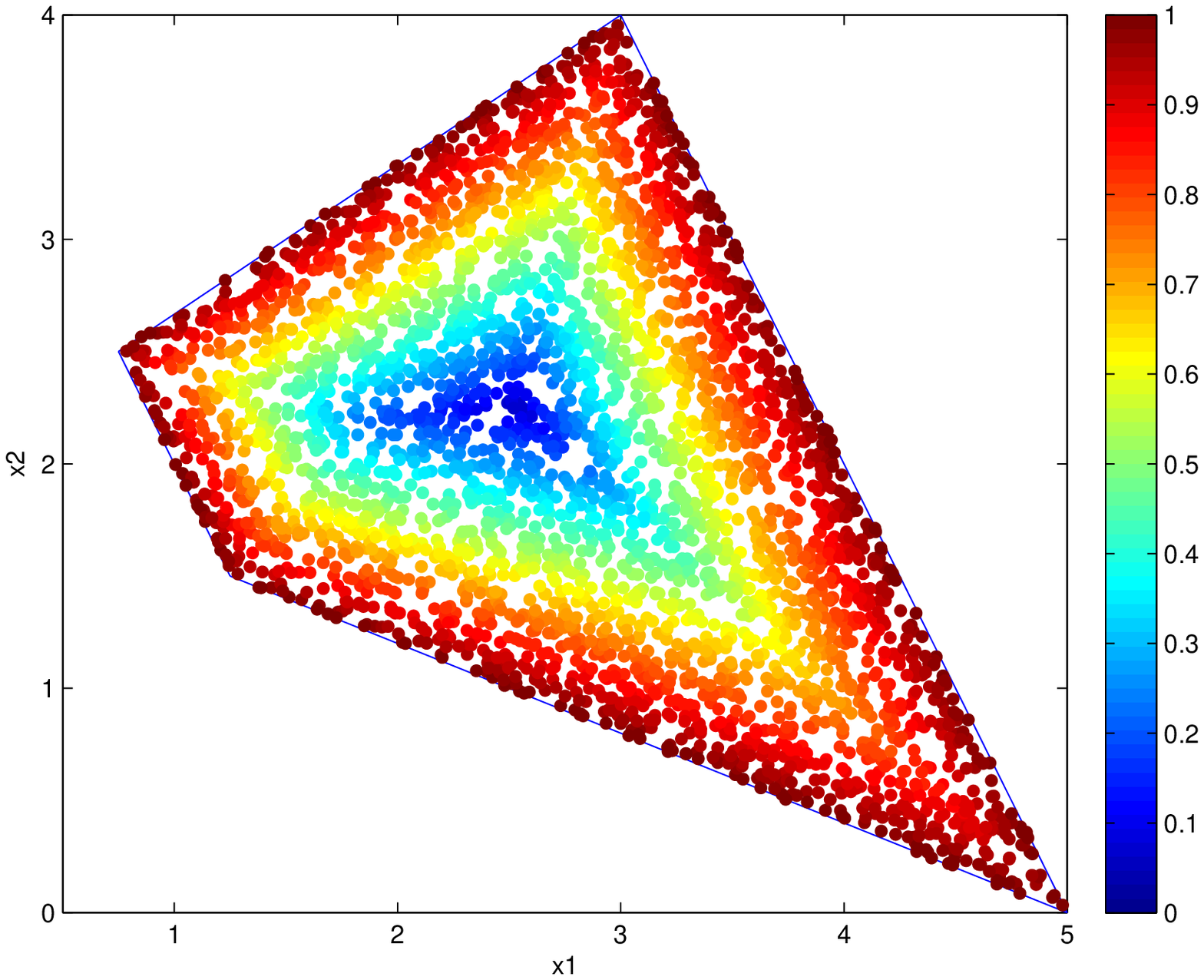}}
\subfigure[Distribution of $\rho_2 - \tilde{\rho}_2$.]
{\label{fig:rho_2_abs_comparison}\includegraphics[width=75mm]{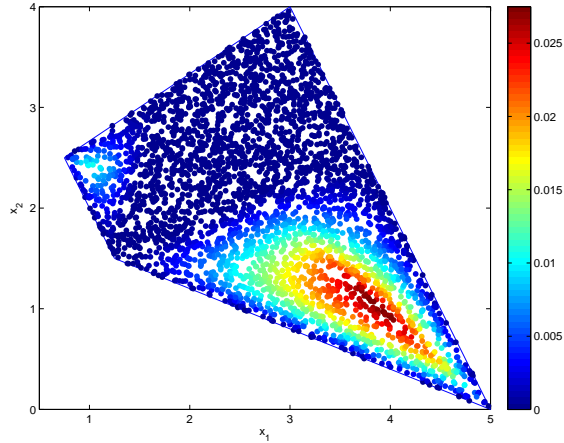}}
\caption{Comparison between $\rho_2$ and $\tilde{\rho}_2$ for Example 1.}
\label{fig:rho_approx}
\end{figure}

Figure~\ref{fig:rho_approx} compares the distribution of $\rho_2$ and $\tilde{\rho}_2$ values inside the polyhedron from Example~\ref{example}.  As shown in Figure~\ref{fig:rho_2_abs_comparison}, the differences are small in magnitude and, by definition, localized to the regions of the polyhedron where the difference between $\bepsilon^i$ and $\tilde\bepsilon^i$ is greatest.  For example, the ``hot'' area to the bottom right of the polyhedron contains those $\hatX$ where the projections $\bepsilon^i$ and $\tilde\bepsilon^i$ with respect to the constraint $2x_1+x_2 \ge 4$ are the most different (also see Figure~\ref{fig:rho_i_illustration}). However, it is possible for the relative difference between $\rho$ and $\tilde\rho$ (in the $p$-norm case) to be arbitrarily large in certain pathological cases, as shown in Example~\ref{ex:pathological}, which has the potential to suggest poor fit when in fact it is good.

\begin{figure}[t]\centering
\subfigure[$\rho=\tilde\rho$.]
{\label{fig:pathological1}\includegraphics[height=6.4cm]{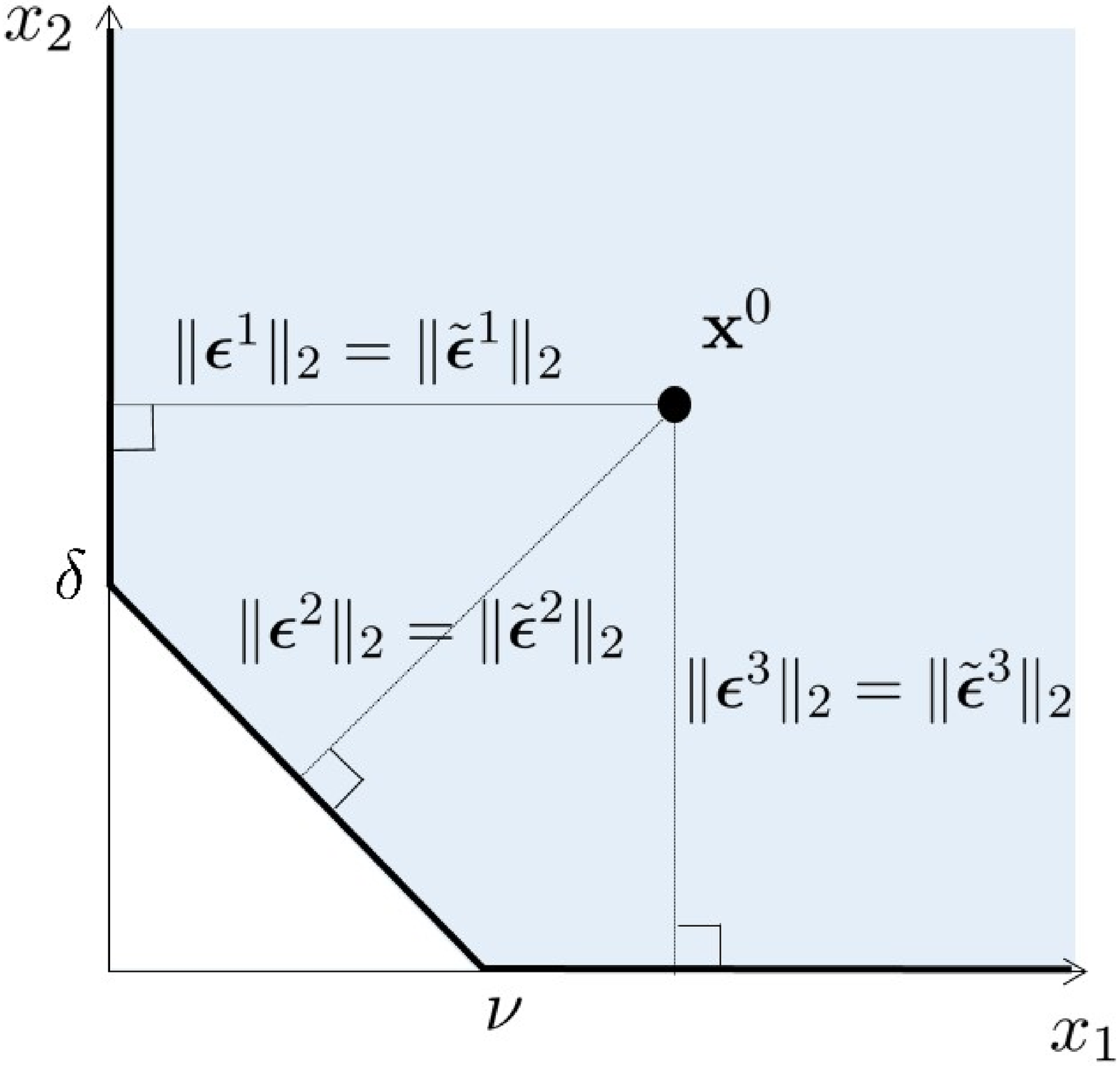}}
\subfigure[$\rho$ and $\tilde\rho$ are different as $\nu$ increases.]
{\label{fig:pathological2}\includegraphics[height=6.6cm, width = 9.7cm]{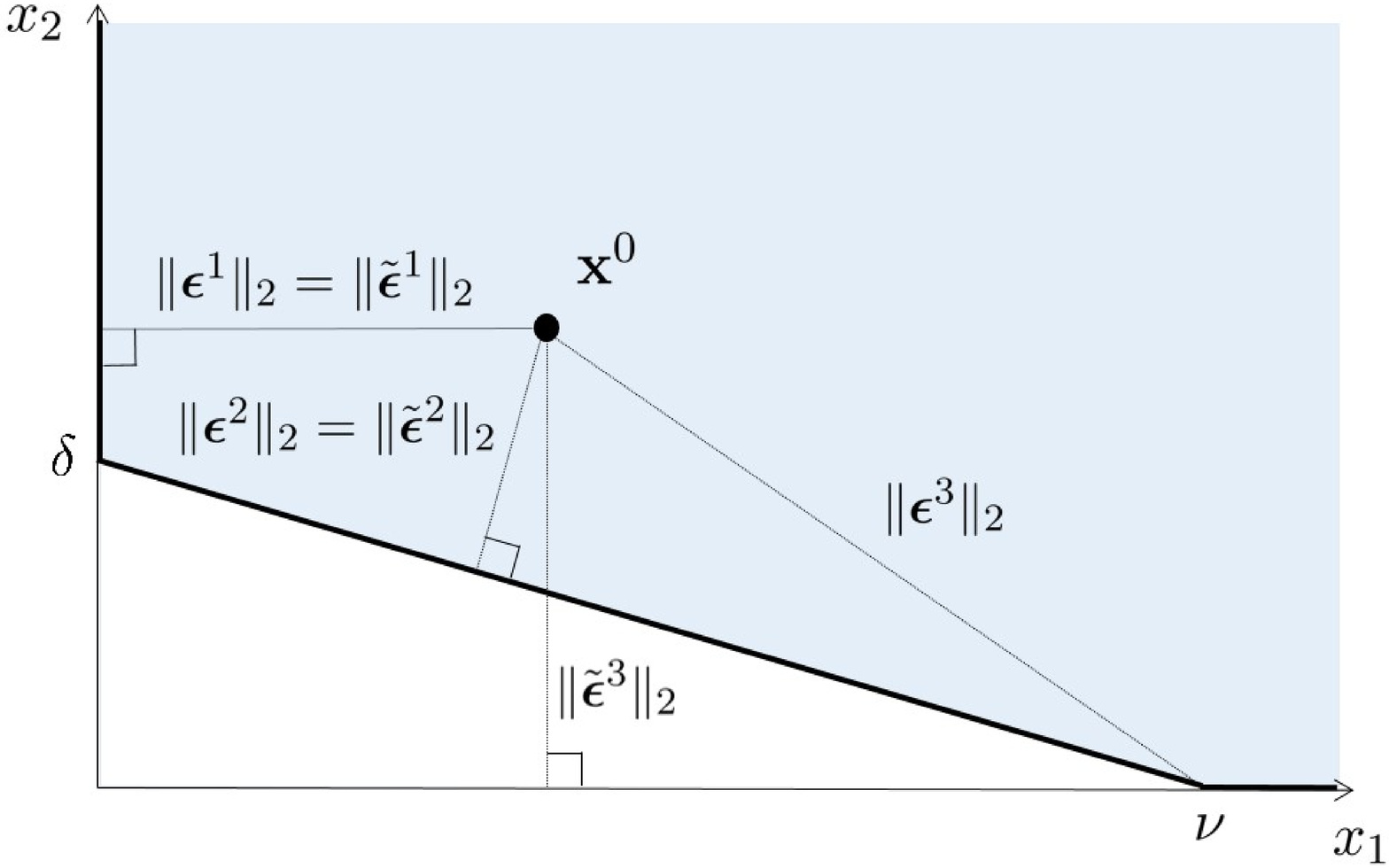}}
\caption{Illustration of Example~\ref{ex:pathological}.}
\label{fig:pathological}
\end{figure}

\begin{example}\label{ex:pathological}
\emph{Let $\bX = \{ (x_1,x_2) \,|\, \frac{1}{\nu} x_1 + \frac{1}{\delta} x_2 \ge 1,\, x_1\ge 0,\, x_2 \ge 0\}$, $\hatX = (1,1)$, $\nu > 0$ and $\delta > 0$. We compare $\rho_2$ and $\tilde\rho_2$. First, consider the case where $\nu = \delta < 1$, shown in Figure~\ref{fig:pathological1}. Since $\bepsilon^i = \tilde\bepsilon^i$ for $i = 1, 2, 3$, $\rho_2 = \tilde\rho_2$. Next, notice that if $\nu > 1 > \delta$, $\|\bepsilon^3\|_2 > \|\tilde\bepsilon^3\|_2$ and thus, $\rho_2 > \tilde\rho_2$ (Figure~\ref{fig:pathological2}).  In fact, consider the extreme case where $\nu \rightarrow\infty$.  Then $\|\bepsilon^1\|_2 = 1$, $\|\bepsilon^2\|_2 \rightarrow 1-\delta$, and $\|\bepsilon^3\|_2 \rightarrow \infty$, and $\rho \rightarrow 1$.  On the other hand, $\|\tilde\bepsilon^1\|_2 = 1$, $\|\tilde\bepsilon^2\|_2 \rightarrow 1-\delta$, and $\|\tilde\bepsilon^3\|_2 = 1$, and therefore $\tilde\rho\rightarrow 1 - \frac{1-\delta}{(1+(1-\delta)+1)/3}  = \frac{2\delta}{3-\delta}$.  Thus, with sufficiently large $\nu$ and sufficiently small $\delta$, we can have $\tilde\rho \approx 0$ while $\rho \approx 1$.}
\end{example}

We expect that for most problems, $\rho$ and $\tilde\rho$ will give similar indications about model-data fit. In addition, this issue can be avoided entirely by using an inverse model that measures error in terms of objective value. Also, when used to support model selection, the difference in $\rho$ values between successive models (i.e., the change in fit) may be a more useful indicator, rather than the absolute value $\rho$. Thus, it may be the case that both $\rho$ and $\tilde\rho$ would serve the purpose equally well. For the $p$-norm case, it may be possible to bound $\tilde\rho$ from below as a function of $\rho$ and geometric properties of the polyhedron $\bX$. Such development is left as future work.

\section{Model Selection and Inverse Optimization}\label{sec:modelSelection}

In fitting a regression model, one is usually faced with the problem of model selection: determining which independent variables results in the best explanatory or predictive performance. Model selection may also look to enforce relationships between certain coefficients (e.g., constrained regression). These issues also exist in inverse optimization. Variable selection in inverse optimization can be thought of as identifying which terms in the objective function are included in the model, and can be formulated as a constrained model where certain cost coefficients $c_j$ are set to 0. In general, a given application may require structural constraints on $\bc$ or $\bepsilon$, which restrict the set of cost vectors that can be returned by the inverse process. Examples of such constraints include fixing certain $c_j$ values, enforcing certain relationships between multiple $c_j$, and constraining $\bepsilon$ so only certain directions of perturbation from $\hatX$ are allowed (e.g., if we want to maintain a certain structure of $\bx^*$ relative to $\hatX$).

In addition to practitioner experience and domain-specific knowledge, the process of model selection in regression is typically guided by goodness of fit (e.g., $R^2$). In inverse optimization, $\rho$ can play a similar role. In the framework that we have developed, $\rho$ can be a tool to help choose the ``best'' model among a set of linear models parameterized by $\bc$, possibly in the presence of structural constraints that enforce prior knowledge on $\bc$ or a particular relationship between $\hatX$ and its projection to the boundary of $\bX$. Below, we discuss several issues regarding model selection in inverse optimization: the impact of structural constraints on $\bc$ and $\bepsilon$ on model estimation, the use and interpretation of $\rho$, and practical considerations regarding the choice of specific $\mathbf{GIO}(\hatX)$ model variant.

\subsection{Model Estimation under Structural Constraints}\label{sec:structuralConstraints}
In the absence of constraints on $\bc$ or $\bepsilon$, Theorem~\ref{thm:GGIO} establishes that $\bc^* = \ba_i / \|\ba_i\|_1$ for some $i$. With the addition of structural constraints of the form $\bc \in \mathcal{C}$ or $\bepsilon \in \mathcal{E}$ to formulation~\eqref{eq:GIO_general}, Theorem~\ref{thm:GGIO} may not be directly applicable. Below, we discuss the effect of such constraints on the applicability of our previously derived closed-form solution and the solution of the resulting inverse problem in general.

\subsubsection{Constraints on $\bc$.}\label{sec:ctsOnC}
When there are restrictions on $\bc$, $\bc^*$ does not in general retain the structure $\ba_i / \|\ba_i\|_1$.  However, the corresponding $\bc$-constrained $\mathbf{GIO}(\hatX)$ problem may still be feasible. For example, an optimal $\bc^*$ may be a strict conic combination of the $\ba_i$ (from the dual feasibility constraints), with the projected point $\hatX - \bepsilon^*$ corresponding to a vertex of the forward polyhedron. As an example, consider a unit square centered at the origin, $\hatX = (0,0)$ and $\bc \in \mathcal{C} = \{\bc\, | \,c_1 = c_2\}$. Clearly $\bc^* \ne \ba_i / \|\ba_i\|_1$ for any $i$, but $\bc = (1/2, 1/2)$ is feasible for the inverse problem. Thus, in general, structural constraints on $\bc$ require the direct solution of the inverse optimization model (i.e., \eqref{eq:GIO_general}), which suggests that the duality gap formulations are preferred in such cases. 

\subsubsection{Constraints on $\bepsilon$.}\label{sec:ctsOnEps}
Constraining $\bepsilon$ affects how $\hatX$ can be projected to the boundary of the forward polyhedron. Since restricting the way $\hatX$ is projected is most relevant in the decision space, we focus on the $p$-norm case. With constraints $\bepsilon\in\mathcal{E}$, $\mathbf{GIO}_p(\hatX)$ can be solved by simply solving
\begin{equation}\label{eq:epsConstraints_p}
\begin{alignedat}{1}
\underset{\bepsilon}{\textrm{minimize}} & \quad \|\bepsilon\|_p \\
\textrm{subject to}
                  & \quad \bA (\hatX - \bepsilon) \ge \bbb, \\
                  & \quad \ba_i' (\hatX - \bepsilon) = b_i, \\
                  & \quad \bepsilon \in \mathcal{E},
\end{alignedat}
\end{equation}
for each $i = 1, \ldots, m$. Formulation~\eqref{eq:epsConstraints_p} determines whether there is a feasible projection to constraint $i$ that satisfies $\bepsilon \in \mathcal{E}$. Let $\hat{I} \subseteq I$ index this set of constraints. As long as $\hat{I}$ is nonempty, the $\bepsilon$-constrained $\mathbf{GIO}_p(\hatX)$ problem is feasible (e.g., $\bc = \ba_i/\|\ba_i\|_1$, $i \in \hat{I}$ is a feasible solution). In fact, solving problem~\eqref{eq:epsConstraints_p} for each $i\in I$ will identify the optimal solution to the $\bepsilon$-constrained $\mathbf{GIO}_p(\hatX)$ problem; it will simply be $\bc^* = \ba_{i}/\|\ba_{i}\|_1$ for $i \in \hat{I}$ with smallest corresponding optimal objective value in~\eqref{eq:epsConstraints_p}. Of course, this solution $\bc^*$ is the closed-form solution from Theorem~\ref{thm:GGIO} applied to the set $\hat{I}$, with corresponding $\by^*$ and $\bepsilon^*$.

\subsubsection{Effect on Calculation of $\rho$.}\label{sec:effect_on_calc}
The addition of structural constraints $\bc \in \mathcal{C}$ or $\bepsilon \in \mathcal{E}$ affects the calculation of $\bepsilon^*$, as outlined in Sections \ref{sec:ctsOnC} and \ref{sec:ctsOnEps}. However, structural constraints may or may not affect the calculation of the denominator of $\rho$, depending on the assumptions underlying the particular application. See Section~\ref{sec:calcOfRho} for discussion on when it is appropriate to adjust or not adjust the denominator calculation, and the associated implications. For now, we discuss the mechanics of making the adjustment, if the modeler chooses to do so.

For the $p$-norm case, taking into account constraints $\bepsilon \in \mathcal{E}$ in the denominator requires no extra effort since finding $\bepsilon^*$ under such constraints already requires finding $\bepsilon^i$ using \eqref{eq:epsConstraints_p} instead of \eqref{eq:eps_p^i}. With constraints on $\bc$ in the absolute and relative duality gap cases, the feasibility of $\epsilon^i_a = (\ba'_i\hatX - b_i)/\|\ba_i\|_1$ and $\epsilon^i_r = \ba_i'\hatX/b_i$ for each $i \in I$ needs to be checked in $\mathbf{GIO}_a(\hatX)$ and $\mathbf{GIO}_r(\hatX)$, respectively; only feasible $\bepsilon^i$ values (i.e., those where the associated $\ba_i$ satisfies the constraints on $\bc$) should be included in the denominator.

\subsection{Guidelines for Practitioners}\label{sec:choiceofGIO}
In this section we discuss additional considerations for the calculation and usage of $\rho$, and provide recommendations regarding the choice of specific $\mathbf{GIO}(\hatX)$ model variant.

\subsubsection{Calculation of $\rho$}\label{sec:calcOfRho}
To determine whether the denominator of $\rho$ should be altered by the structural constraints, the practitioner needs to determine whether these constraints change the set of candidate cost vectors for the underlying problem. Structural constraints should be applied in the calculation of the denominator of $\rho$ only when the set of candidate cost vectors is reduced as a result of the constraints. For example, in our numerical experiments (Section~\ref{sec:experimentalResults}), negative cost vectors are not physically meaningful and thus should be eliminated from consideration in the denominator. On the other hand, in Example 2, the $\bc$ vector is unconstrained but we may want to evaluate the quality of fit provided by the best non-negative $\bc$ vector; in this case, the denominator should not be affected by the non-negativity constraint. 

Note that if the denominator calculation is not modified in the presence of structural constraints, the non-negativity property of $\rho$ may not hold anymore as $\|\bepsilon^*\|$ may be larger than $\frac{1}{m} \sum_{i=1}^{m} \|\bepsilon^i\|$. That is, the fit of any of the $m$ possible vectors in the unconstrained case is better than that of the optimal constrained one. The analogous result in regression is that $R^2$ may be negative when used to evaluate a constrained regression model. There is also an analogous situation in regression to when the the $\rho$ calculation should be adjusted due to structural constraints. Consider regression through the origin, where only in models with a zero intercept are considered. In this case, it is recommended that the calculation of $R^2$ be modified to use a denominator of $\sum_{q=1}^{Q} y_q^2$ instead of $\sum_{q=1}^{Q} (y_q-\bar{y})^2$, effectively imposing a constraint $\beta_0 = 0$ in the denominator \citep{Sprecher94}.

\subsubsection{Usage of $\rho$}\label{sec:usageOfRho}
Analogous to $R^2$ in regression, there is no universal threshold that determines acceptable fit when using $\rho$. In addition to considering its value, the user should also consider the change in $\rho$ when comparing different models, with the general caveat that overfitting should be balanced against obtaining a high value of $\rho$.

\subsubsection{Choice of $\mathbf{GIO}(\hatX)$ variant}
The choice of a specific $\mathbf{GIO}(\hatX)$ model to use may be dictated by the application. For instance, if the modeler aims to recover a solution $\bx^*$ that is optimal for the inversely optimized $\bc^*$ and ``close'' to the initial $\hatX$, a $p$-norm model would be most natural. If the goal is to recover the objective function value achieved by $\hatX$, then a duality gap model would be a more appropriate choice.

When there are constraints on $\bc$ that prevent the modeler from applying the closed-form solution in Theorem~\ref{thm:GGIO}, our general recommendation is to use one of the duality gap models since they can be solved as linear programs under mild assumptions that are often
satisfied in real applications (e.g., non-negative costs). When there are constraints on $\bepsilon$, on the other hand, $\mathbf{GIO}_p(\hatX)$ is naturally better suited; in this case, the structure of the solution is characterized by Theorem~\ref{thm:GGIO} applied to the set $\hat{I}$ of feasible constraints as discussed in Section~\ref{sec:ctsOnEps}.

\section{Numerical Studies}\label{sec:experimentalResults}
In this section, we demonstrate the use of $\rho$ in model estimation and evaluation for two different applications that have previously been studied in the inverse optimization literature: production planning \citep{Troutt06} and cancer therapy \citep{Chan14}.

\subsection{Production Planning}\label{sec:APP}
\citet{Troutt06, Troutt08} posit that the application of mathematical programming for production planning in practice is complicated by the lack of appropriate parameter estimates. However, given past production data, inverse optimization may provide a way to overcome this difficulty. We demonstrate the application of our inverse optimization framework to the aggregate production planning (APP) problem from \citet{Troutt06}. In particular, we show how initial managerial estimates for certain cost parameters can be updated through an iterative process of model refinement and evaluation using $\rho$. A more accurate APP model can aid prospective decision-making in the production planning process.

\subsubsection{Problem Description}

The APP model (forward problem) aims to find the optimal production quantities for each quarter $h$ of a year to minimize costs and meet demand. In this model, there are five decision variables: regular time production, overtime, idle time, inventory, and backorders, all expressed in terms of direct labour hours. For example, the number of items in inventory is expressed in terms of the number of labour hours needed to produce them. The APP model is formulated as the following linear program:
\begin{subequations}\label{eq:APP}
\begin{align}
\quad \underset{\bx}{\text{minimize}} & \quad \; \sum_{h=1}^{4} \sum_{j = 1}^{5} c_j x_{h, j} \\
\text{subject to} & \quad x_{h-1, 4} - x_{h-1, 5} - x_{h, 4} + x_{h, 5} + x_{h, 1} + x_{h, 2} = D_{h}, \;\; h = 1, \dots, 4, \label{eqn:inventoryBalancing}\\
& \quad  x_{h, 1} + x_{h, 3} = A_1, \quad h = 1, \dots, 4, \label{eqn:regularTimeMax} \\
& \quad x_{h, 2} \le A_2, \quad h = 1, \dots, 4, \label{eqn:overtimeMax} \\
& \quad x_{h, j} \ge 0, \quad h = 1, \dots, 4, \; j = 1, \dots, 5, \label{eqn:APP_nonNeg}
\end{align}
\end{subequations}
where the decision variables $x_{h, 1}, x_{h, 2}, x_{h, 3}, x_{h, 4}, x_{h, 5} \; (h = 1, \ldots, 4)$ represent the per-quarter regular-time production hours, overtime production hours, idle-time hours, hours of production stored in inventory and hours of production backordered, respectively, and $c_1, c_2, c_3, c_4, c_5$ are the corresponding costs per hour. The model parameters are $D_h$, the forecasted demand in hours; $A_1$, the maximum number of regular-time hours available per quarter; and $A_2$, the maximum number of overtime hours available per quarter. We define the parameters $x_{0, 4}$ and  $x_{0, 5}$ as the number of hours left in inventory and backordered, respectively, from the previous year. 

We aim to obtain per-hour cost estimates that allow the model to replicate the total production cost associated with a given solution as closely as possible. Thus, we choose an inverse model based on minimizing the error in objective value, in particular the absolute duality gap formulation~\eqref{eq:ailop_a}; the complete model is shown in Appendix~\ref{sec:APP_inverseModel}.

\subsubsection{Data}
We use one year of past demand data from \cite{Troutt06}'s Table 1 as the \emph{realized past demand}, $\bD$.\footnote{In contrast, Troutt et al. use these values as \emph{forecasts}, rather than demand realizations.} To generate $\hatX$, we solve formulation~\eqref{eq:APP} using $\bD$, assuming costs $c_1 = \$14, c_2 = \$21, c_3 = \$8, c_4 = \$4, c_5 = \$17$, and then perturb the resulting optimal solution similarly to \citet{Troutt06} (see Appendix~\ref{sec:APP_observation} for details). 

\subsubsection{Evaluating cost assumptions using $\rho_a$}
We illustrate how a manager may use $\rho$ to guide an iterative process of model and assumption refinement. A series of four models is presented below, each with different assumptions on $\bc$ and representing a refinement of a previous model. A summary of the results is given in Table \ref{tab:APP_results}; additional details regarding the implementation of assumptions on $\bc$ are given in  Appendix~\ref{sec:appendix_APP_scaling}. Because we employ the absolute duality gap model, all models are evaluated using $\rho_a$~\eqref{eq:singleObsRhoGeneralGIO_a}. We calculate the terms $\epsilon^i_a$ in the denominator of $\rho_a$ using equation \eqref{eq:eps_i_a_def}. Because we know that all relevant solutions should have non-negative costs, as recommended in Section~\ref{sec:usageOfRho} we choose to compute both the numerator and the denominator of $\rho$ with the constraints $c_j \ge 0.0001$ $\forall j$. Thus, we check the feasibility of each $\epsilon^i_a$ by solving the inverse model with the constraint $\epsilon_a = \epsilon^i_a$ and include only the feasible $\epsilon^i_a$ values in the denominator of $\rho_a$.

In Table \ref{tab:APP_results}, Model 1 is the most constrained, representing the manager's prior belief about the costs; the corresponding $\rho_a$ of 0.426 suggests that this belief does not match the observed solution, particularly because the estimated cost of production is too low and the estimated cost of idle time is too high for an observed solution where the regular-time production is not at capacity. Model 2 relaxes one of the manager's prior constraints; the substantial improvement in fit combined with a low estimate for $c_3$ suggests that even though the manager believed idle time was expensive ($c_3 \ge \$24$ in Model 1), production decisions were made as if idle time was very cheap ($c_3 = \$0.0035$). Adjusting the constraints to allow the cost of regular production to be higher in Model 3, we obtain an even better $\rho_a$. An upper bound on the goodness of fit for the given solution is obtained by considering a completely unconstrained model, Model 4. This model results in almost perfect fit, but the estimated cost vector includes an estimate for the inventory cost that is several orders of magnitude larger than the other costs. The manager may therefore decide that Model 3 provides a sufficiently good cost estimate for prospective production planning.

\begin{table}[ht]
\centering
\begin{tabular}{cccr} \hline
Model & \multicolumn{1}{c}{Constraints on $\bc$} & $\bc^*$ & \multicolumn{1}{c}{$\rho_a$} \\ \cmidrule(l){1-1}\cmidrule(l){2-2}\cmidrule(l){3-3}\cmidrule(l){4-4}
1 & $2c_2 = 21c_4$, $c_1 \le 3c_4$, $c_3 \ge 12c_4$ & $[6, 21, 24, 2, 25]$ & $\mathbf{0.426}$ \\
2 & $2c_2 = 21c_4$, $c_1 \le 3c_4$ & $[6, 21, 0.0035, 2, 6]$ & $\mathbf{0.846}$ \\
3 & $4c_2 = 21c_4$, $c_1 \le 3c_4$  & $[12, 21, 1.5, 4, 10.5]$ & $\mathbf{0.906}$ \\
4 & None & $[42, 21, 21, 209895, 21]$ & $\mathbf{0.999}$ \\
\hline
\end{tabular}
\caption{Cost estimates from different inverse APP models, scaled to $c_2 = 21$, and the corresponding $\rho_a$ values. (See Appendix~\ref{sec:appendix_APP_scaling} for the unscaled $\bc$ vector.)}
\label{tab:APP_results}
\end{table}

The substantial fit improvement between Models 1 and 3 suggests that the manager's initial cost estimates did not represent well how the observed production plan was generated. The remaining (unexplained) model-data error may be attributed to several factors including sub-optimal historical production planning, additional inaccurate cost assumptions, and inaccurate past demand forecasts. These factors can easily be explored further by the manager in a manner similar to the one demonstrated. In Model 4, even though the high $\rho_a$ value suggests excellent fit, the estimated $\bc$ is far (in terms of norm) from the true $\bc$ since our inverse model is focused on minimizing error in the objective value space (i.e., duality gap). However, constraints of the form used in Models 1 to 3 can help to encourage inverse solutions $\bc$ that are close to a given prior cost vector.

\subsection{Cancer Therapy}\label{sec:cancertherapy}
Intensity-modulated radiation therapy (IMRT) is a cancer treatment technique that uses high energy x-ray beams to deliver radiation to a tumor. IMRT treatments can be designed using an optimization model with a composite objective function comprising a weighted sum of several organ-specific objectives that trade off between tumor dose and healthy organ dose. A major challenge in IMRT treatment plan optimization is determining the right objective functions and corresponding weights necessary to generate a clinical quality treatment. For example, the relative improvement that may be gained through the addition of a particular objective is not known \emph{a priori}.  In practice, trial-and-error is used.  Recently, \citet{Chan14} demonstrated that appropriate weights could be determined from past clinically approved treatments of prostate cancer using inverse optimization. These weights could then be used to train statistical models that predict weights from patient anatomy to create treatments for \emph{de novo} patients~\citep{Lee13,Bout15}. However, the overall process is heavily dependent on the initial choice of objective functions.
In this section, we demonstrate the use of $\rho$ to measure which objectives contribute the most to re-creating a clinical IMRT treatment plan. Using a well-known treatment evaluation tool, we also provide a context-specific validation of $\rho$ as a measure of the goodness of fit. 

\subsubsection{Problem description}
The aim of the IMRT optimization problem in prostate cancer is to deliver a sufficiently high dose of radiation to the prostate while minimizing dose to the neighboring organs-at-risk (OARs). The primary disease site (prostate) is referred to as the clinical target volume (CTV), and typical OARs include the bladder, rectum, and left and right femoral heads.

We employ the linear forward and inverse multi-objective optimization formulations of~\cite{Chan14}, which are included in Appendix~\ref{appendix:IMRT} for completeness. The forward problem has an objective function of the form $\balpha'\bC\bx$, where $\bC$ is a matrix of (row-wise) objectives and $\balpha$ is the vector of objective weights to be estimated. In other words, this model is equivalent to a linear program with the cost vector $\bc$ constrained to be a conic combination of the rows of $\bC$. As such, there is no guarantee that any of the vectors $\ba_i$ (appropriately scaled) can be written as $\bC'\balpha$ and thus the inverse problem needs to be solved directly. The specific inverse model employed is an instance of the relative duality gap formulation, so goodness of fit is measured by $\rho_r$.

\subsubsection{Data}
We used a clinical prostate cancer treatment plan from Princess Margaret Cancer Centre, i.e., a particular patient's radiation dose distribution, as input. The beamlet intensities (i.e., the observed solution $\hatX$) were not available, but having the complete dose distribution is all that is needed to compute the vector of objective function values $\bC\hatX$, which is the input required for a multi-objective inverse optimization problem.

\subsubsection{Evaluating objective functions using $\rho_r$}
We illustrate how a treatment planner may use $\rho_r$ to evaluate how different objectives 
improve the ability of the model to re-create the given dose distribution.  Following~\citet{Chan14}, two types of objective functions are considered: (a) a piece-wise linear function that penalizes dose to any part of the OAR above a certain dose threshold, and (b) a maximum dose function. We consider objectives of type (a) with thresholds of $0, 10, \ldots, 70$ Gray (Gy) for the bladder and rectum, and objectives of type (b) for the femoral heads.

Unlike the production planning application in the previous section, we test different subsets of the 18 objectives by starting with a simple set and adding more objectives to increase $\rho_r$, which emulates the current practice of treatment planning. Note that specifying the subset of the objectives to use in the inverse formulation is equivalent to fixing weight values for the other objectives to zero. However, such constraints do not affect the calculation of the denominator of $\rho_r$ because for any subset of objectives we want to measure its fit relative to the complete 18-objective model (see Section~\ref{sec:effect_on_calc}). We calculate the terms $\epsilon^i_r$ in the denominator of $\rho_r$ using equation \eqref{eq:eps_i_r_def}. The purpose of this experiment is not to provide a prescriptive scheme for determining which objectives should be added or removed. Rather, the goal is to observe how $\rho_r$ varies as the treatment planner explores the explanatory power of different objectives used in concert or alone.

\begin{table}\setlength{\tabcolsep}{2.5pt}
\begin{center}
\resizebox{\columnwidth}{!}{\begin{tabular}{lccccccccccccccccccc}\hline
\multirow{ 2}{*}{Model}  & \multicolumn{2}{c}{$\theta=0$} & \multicolumn{2}{c}{$\theta=10$}&\multicolumn{2}{c}{$\theta=20$}&\multicolumn{2}{c}{$\theta=30$}&\multicolumn{2}{c}{$\theta=40$}&\multicolumn{2}{c}{$\theta=50$} &\multicolumn{2}{c}{$\theta=60$}&\multicolumn{2}{c}{$\theta=70$}& \multicolumn{2}{c}{Max} & \multirow{ 2}{*}{$\rho_r$} \\\cmidrule(l){2-3}\cmidrule(l){4-5}\cmidrule(l){6-7}\cmidrule(l){8-9}\cmidrule(l){10-11}\cmidrule(l){12-13}\cmidrule(l){14-15}\cmidrule(l){16-17}\cmidrule(l){18-19}
	     & B & R & B & R & B & R & B & R& B & R & B & R & B & R & B & R & LF & RF &  \\\hline
1 &  &  &  &  &  &  &  &  &  &  &  &  &  &  & \cellcolor{gray!25}{0.667} & \cellcolor{gray!25}{0.333} &  &  & \textbf{0.724} \\
2 &  &  &  &  &  &  &  & &  &  &  &  &  &  & \cellcolor{gray!25}{0.604} & \cellcolor{gray!25}{0.394} & \cellcolor{gray!25}{0.001} & \cellcolor{gray!25}{0.001} & \textbf{0.737} \\
3 & \cellcolor{gray!25}{0.591} & \cellcolor{gray!25}{0.063} &  &  &  &  &  &  &  &  &  &  &  &  & \cellcolor{gray!25}{0.326} & \cellcolor{gray!25}{0.001} & \cellcolor{gray!25}{0.010} & \cellcolor{gray!25}{0.009} & \textbf{0.947} \\
4
& \cellcolor{gray!25}{0.582} & \cellcolor{gray!25}{0.058}
& \cellcolor{gray!25}{\hskip 0.15in - \hskip 0.05in} & \cellcolor{gray!25}{\hskip 0.15in - \hskip 0.05in}
& \cellcolor{gray!25}{\hskip 0.15in - \hskip 0.05in} & \cellcolor{gray!25}{\hskip 0.15in - \hskip 0.05in}
& \cellcolor{gray!25}{\hskip 0.15in - \hskip 0.05in} & \cellcolor{gray!25}{\hskip 0.15in - \hskip 0.05in}
& \cellcolor{gray!25}{\hskip 0.15in - \hskip 0.05in} & \cellcolor{gray!25}{\hskip 0.15in - \hskip 0.05in}
& \cellcolor{gray!25}{0.059} & \cellcolor{gray!25}{0.008}
& \cellcolor{gray!25}{\hskip 0.15in - \hskip 0.05in} & \cellcolor{gray!25}{\hskip 0.15in - \hskip 0.05in}
& \cellcolor{gray!25}{0.275} & \cellcolor{gray!25}{\hskip 0.15in - \hskip 0.05in}
& \cellcolor{gray!25}{0.009} & \cellcolor{gray!25}{0.009}
& \textbf{0.948} \\\hline
\multicolumn{20}{l}{\footnotesize B = bladder, Re= rectum, LF = left femoral head, RF = right femoral head.   $\theta$ represents the dose threshold in objective type (a). }\\
\multicolumn{20}{l}{\footnotesize Gray cells indicate the objectives included for each model. Dashes indicate a weight of zero.}
\end{tabular}}
\end{center}
\caption{Results from inverse models with different sets of objectives.}
\label{tab:imrt}
\end{table}

Table~\ref{tab:imrt} shows the weights and the $\rho_r$ values for models with different choices of objectives. Gray cells indicate which objectives were included in the model while the values of the cells indicate the corresponding objective function weight. Model 1 is the simplest, with two objectives that penalize dose above 70 Gy to the bladder and rectum; the corresponding $\rho_r$ is $0.724$. Model 2 includes two additional objectives for the femoral heads, but the fit improves only slightly. Adding the mean dose objectives (equivalent to dose threshold of 0 Gy) for the bladder and rectum in Model 3, we see $\rho_r$ increases to $0.947$, suggesting a good fit with the historical treatment plan. Finally, Model 4 considers all of the 18 objectives. For this model, the value of $\rho_r$ is almost the same as for Model 3 (0.948). Interestingly, the corresponding weight vector has only seven nonzero components. Overall, for this specific patient, the bladder seems to be the most important organ to spare.

\begin{figure}\centering
\subfigure[Bladder DVH.]
{\label{fig:DVH_Blad}\includegraphics[height=70mm]{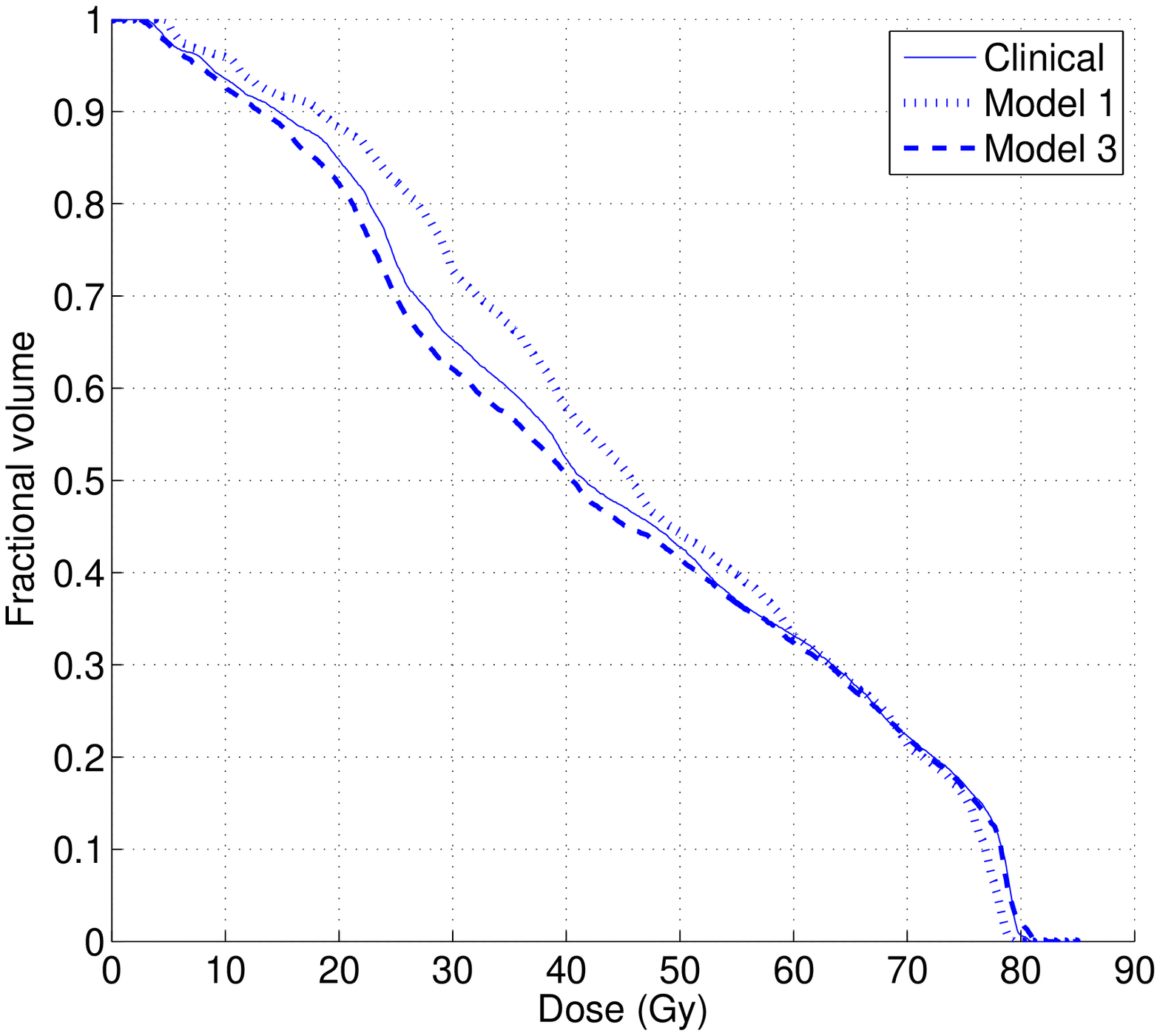}}
\subfigure[Rectum DVH.]
{\label{fig:DVH_Rectum}\includegraphics[height=70mm]{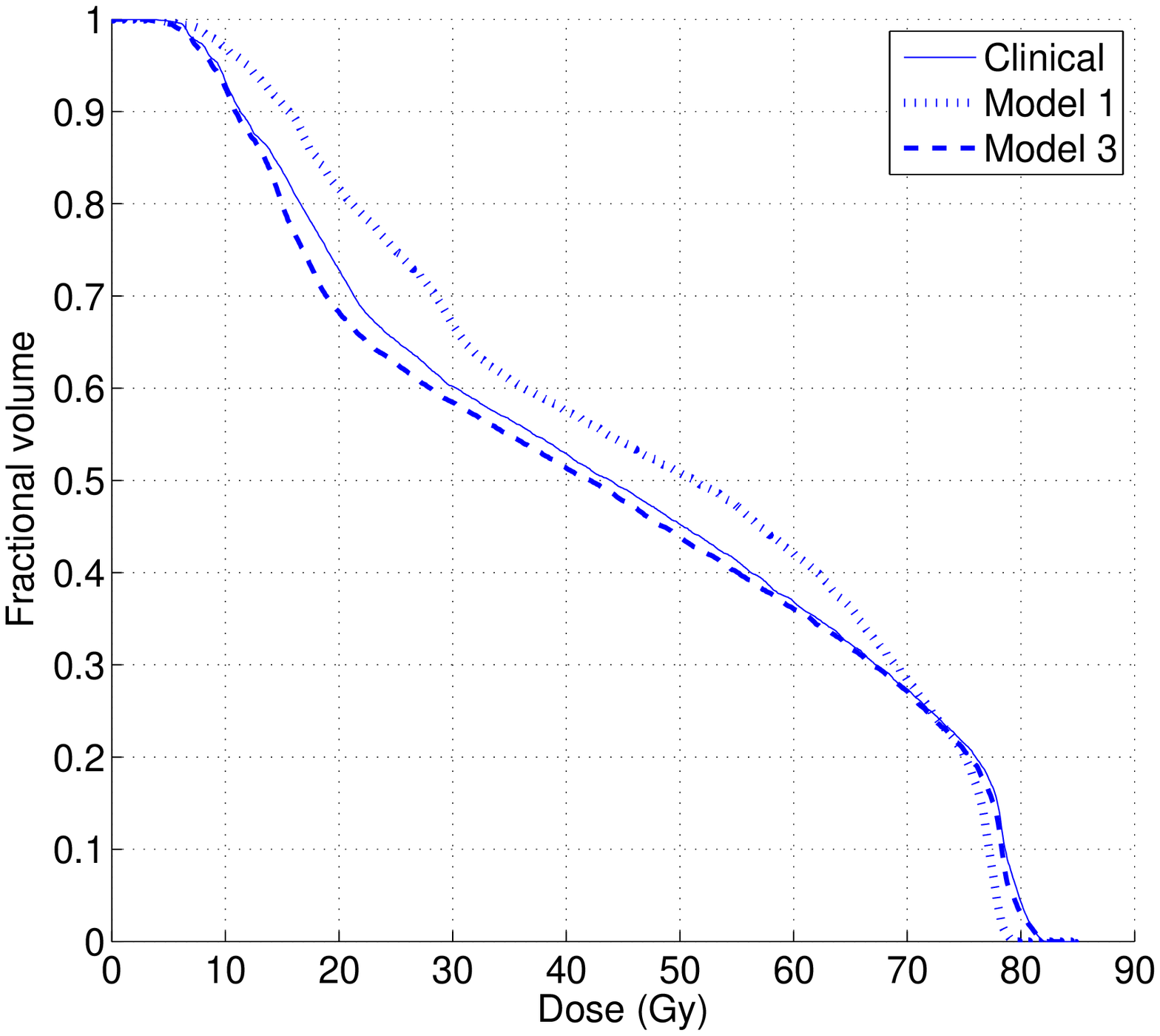}}
\caption{DVHs from the clinical, Model 1, and Model 3 plans.}
\label{fig:rho_DVH}
\end{figure}

To provide a context-specific validation of $\rho_r$ as a measure of goodness of fit, we solve the forward model using the inversely optimized weight vector $\balpha^*$ and compare the resulting dose distribution with the clinical dose distribution using a dose-volume histogram (DVH). A DVH is the clinical standard for visualizing and evaluating treatment quality. For each organ, it shows what fraction of the organ receives a certain dose or higher.  Figure~\ref{fig:rho_DVH} shows the bladder and rectum DVHs corresponding to the clinical plan, the Model 1 (two objectives) plan, and the Model 3 (six objectives) plan. For both organs, the DVH from the six-objective model more closely approximates the clinical DVH than the DVH from the two-objective model does. These results suggest that $\rho_r$ provides an indication of how representative a dose distribution generated with a particular set of objectives and weights is of a clinical dose distribution.

\subsubsection{Sparsity and variable selection}
The process of adding objectives and measuring the resulting goodness of fit in the previous section can be seen as analogous to variable selection in regression. Recall that the inverse solution to the 18-objective model (Model 4) has only seven nonzero objective function weights, indicating that not all objectives were needed to fit the observed treatment plan. Furthermore, while our inverse approach does not explicitly focus on sparsity, it nevertheless seems to generate a sparse solution.  Motivated by this observation, we investigate whether it is possible to find an even sparser solution than $\balpha^*$ while minimally degrading $\rho_r$. To do so, we add a regularization term $\|\balpha\|_1$ to the original objective function of the inverse formulation in Appendix~\ref{appendix:IMRT}, with a weight of $\lambda \in [0,1]$ for the regularization term and a weight of $1-\lambda$ for the original objective, analogous to lasso for regression \citep{Tibshirani94}. The number of nonzero components of $\balpha^*$ drops to six at $\lambda = 0.04$ ($\rho_r = 0.948$) and to five at $\lambda = 0.07$ ($\rho_r = 0.947$). For $\lambda = 0.1$, $\balpha^*$ has only four nonzero values, corresponding to the bladder and rectum mean dose objectives, and the left and right femoral head objectives. The $\rho_r$ value is $0.947$, almost identical to that of Model 3. The same solution is optimal for $\lambda$ values from 0.1 to 0.9.

Using a regularization approach may result in a potentially sparser inverse solution with minimal degradation in goodness of fit. However, evaluating a few other patients using the same process, we find that the level of sparsity of the non-regularized solution varies quite a bit. Interestingly, for some patients, the non-regularized model generates a sparse solution and sparsity actually worsens when including the regularization objective with a nonzero $\lambda$. We believe this result is due to the fact that, as noted in Section \ref{sec:relative}, the relative duality gap model under the assumption of non-negative weights (cost vector), has an objective that is already a weighted $\ell_1$ regularization of the weight vector. Thus, creating an overall objective function that is a convex combination of the original objective function and $\|\balpha\|_1$ may actually generate less sparse solutions than having the original objective alone. 

\section{Conclusions and Future Work}\label{sec:conclusions}
In this paper, we introduce a framework for imputing optimization model parameters from imperfect data using inverse linear optimization, integrating parameter estimation to minimize error in the data-model fit and the measurement of the resulting goodness of fit. Our general inverse optimization model has a closed-form solution and tractable model variants that include previous models from the literature. Our goodness-of-fit metric, the \emph{coefficient of complementarity}, has similar properties as $R^2$ from regression, is easy to interpret, and is broadly applicable. It is the first context-free, general goodness-of-fit metric for inverse optimization and has the potential to support broader uptake of inverse optimization as a tool for optimization model estimation and selection. The choice of inverse optimization model variant should depend on the application setting, with the duality gap models being better suited in the case of constraints on $\bc$ and the $p$-norm models being better suited in the case of constraints on $\bepsilon$.

A natural next step for the work presented in this paper is to extend to the simultaneous treatment of multiple data points, including infeasible points, which would enable broader uptake and applicability of our framework. As our framework draws numerous parallels with regression (e.g., existence of closed-form optimal solutions, goodness-of-fit metric with similar properties), the case of multiple points offers a rich setting for further adaptation of regression concepts. We also hypothesize that our framework can be extended to conic problems, but will likely require new theoretical development to deal with the projection to the complement of a general convex set and the corresponding challenges in computing $\rho$. Finally, this paper takes a geometric perspective of goodness of fit in inverse optimization. The development of a statistical foundation is a worthwhile topic for future study.

\section*{Appendix}
\renewcommand{\thesubsection}{\Alph{subsection}}

\subsection{Aggregate Production Planning (APP) Experimental Details}

\subsubsection{APP Inverse Model}\label{sec:APP_inverseModel}

\begin{subequations}\label{eq:inverse_APP}
\begin{align}
\quad \underset{\bc, \epsilon_a, \bgamma, \blambda, \bmu}{\text{minimize}} & \quad \epsilon_a \\
\text{subject to} & \quad  -\gamma_h + \gamma_{h+1} \le c_4, \quad h = 1, \dots, 3, \label{eqn:APP_1}\\
& \quad  -\gamma_{4} \le c_4, \\
& \quad  \gamma_{h} - \gamma_{h+1} \le c_5, \quad h = 1, \dots, 3, \\
&  \quad \gamma_{4} \le c_5, \\
& \quad  \gamma_h + \lambda_h \le c_1, \quad h = 1, \dots, 4, \\
& \quad  \gamma_h - \mu_h \le c_2, \quad h = 1, \dots, 4, \\
& \quad  \lambda_h \le c_3, \quad h = 1, \dots, 4, \label{eqn:APP_7}\\
& \quad  \mu_h \ge 0, \quad h = 1, \dots, 4, \label{eqn:APP_7} \\
& \quad  \sum_{h=1}^{4} \sum_{j=1}^{5} (c_j x^0_{h, j}) = (D_1 - x_{0,4} + x_{0,5}) \gamma_{1} + \sum_{h=2}^4 D_h\gamma_h \nonumber \\
& \quad \quad  \quad \quad  \quad \quad \quad \;   +  \sum_{h=1}^4 (A_1\lambda_h- A_2\mu_{h}) + \epsilon_a, \label{eqn:strongDualityAPP} \\
& \quad \sum_{j=1}^{5} c_j = 1, \label{eqn:normalization}\\
& \quad  c_j \ge 0, \quad j = 1, \dots, 5. \label{eqn:nonZero}
\end{align}
\end{subequations}

The inputs to the inverse model~\eqref{eq:inverse_APP} are the observed production plan ($\hatX$), the number of hours in inventory/backordered from last year ($x_{0,4}, x_{0,5}$), and the demand realization ($\bD$). The variables $\bgamma$, $\blambda$, $\bmu$ are the vectors of dual variables corresponding to constraints~\eqref{eqn:inventoryBalancing}, \eqref{eqn:regularTimeMax} and \eqref{eqn:overtimeMax}, respectively. Equations~\eqref{eqn:APP_1}--\eqref{eqn:APP_7} represent dual feasibility, and equation~\eqref{eqn:strongDualityAPP} is the relaxed strong duality constraint. We assume that all costs are non-negative, allowing normalization to be enforced by~\eqref{eqn:normalization} and \eqref{eqn:nonZero}.

\subsubsection{Given APP Solution}\label{sec:APP_observation}
Table~\ref{tab:observedSolution_APP} lists the data used in the APP example from Section~\ref{sec:APP}. Columns 2--6 correspond to $\bx^0$ and column 7 shows the given demand realization $D_h$, $h = 1, \dots, 4$. The numbers in the brackets show the original, non-perturbed optimal solution to the forward problem with demand $\bD$. Note that the number of hours left in inventory from the previous year ($x_{0,4}$) and the number of hours backordered from the previous year ($x_{0, 5}$) are different for the original and perturbed solution. That is, $x_{0,4} = 0$ and $x_{0, 5} = 0$ for the original non-perturbed solution, while $x_{0,4}$ is 0 and $x_{0, 5} = 1,183.8$ for the perturbed solution.

The perturbed solution is produced as follows. We first perturb the regular time production by adding either a random number $U[0, 0.02]$ times the idle time (if original production time is 0) or a random number $U[-0.03, 0.03]$ times the original regular time production, and then adjusting all the remaining quantities to ensure a feasible solution. Idle time is adjusted to reflect the increase or decrease in regular time hours, ensuring feasibility for constraint~\eqref{eqn:regularTimeMax}. We do not adjust overtime hours, and so constraint~\eqref{eqn:overtimeMax} is automatically satisfied.
If the sum of demand and the initial backorder quantity exceeds the sum of the initial inventory, the perturbed production level and the overtime, then the backorder value is increased by the corresponding amount, while the perturbed inventory level is set to 0. Otherwise, the inventory is increased and the backorder is set to 0. Thus, we ensure feasibility of the perturbed solution for constraint~\eqref{eqn:inventoryBalancing}.

\begin{table}[ht]
\centering
\small
\begin{tabular}{crrrrrc} \hline
Quarter & \multicolumn{5}{c}{Production Plan} & Demand \\ \cmidrule(l){1-1}\cmidrule(l){2-6}\cmidrule(l){7-7}
 & \multicolumn{1}{c}{Regular-time} & \multicolumn{1}{c}{Over-time} & \multicolumn{1}{c}{Idle-time} & \multicolumn{1}{c}{Inventory} & \multicolumn{1}{c}{Backorder}  \\ \cmidrule(l){2-6}
1 & 34,430.9 (35000) & 3,500 (3500) & 569.1 (0) & 0 (0) & 13,453.0  (11700) & 50,200  \\
2 & 34,114.2 (35000) & 3,500 (3500)  & 885.8 (0) & 0 (0) & 19,763.8 (17125) & 43,925 \\
3 & 29,056.1 (29675) & 0 (0) & 5,943.9 (5325) & 0 (0) & 3,257.8 (0) & 12,550\\
4 & 21,112.8 (20708) & 0 (0) & 13,887.2 (14292) & 0 (0) & 2,852.9 (0) & 20,708\\ \hline
\end{tabular}
\caption{The observed APP solution from Section 6.1. 
The numbers in brackets show the non-perturbed solution. All units in hours.}
\label{tab:observedSolution_APP}
\end{table}

\subsubsection{Scaling}\label{sec:appendix_APP_scaling}
Assumptions on $\bc$ are implemented via the rescaling method from \citet{Troutt06}. That is, assuming $J_0 \subset \{1, 2, 3, 4, 5\}$ is a set indexing those costs that are assumed to be known a-priori ($\bar{c}_j$), we add to the inverse model the linearized versions of the constraints $c_l/c_j = \bar{c}_l/\bar{c}_j$ for $l, j \in J_0$, together with $c_j \ge 0.0001$ $\forall j$ to prevent division-by-zero errors. We scale all estimated cost vectors so that $c_2 = 21$, assuming the manager is confident in the estimate of this cost a-priori. Table \ref{tab:APP_results_nonscaled} shows the solutions to Model 1 to 4 without the scaling. Note that the choice of scaling can impact the interpretation of the results and thus should be chosen carefully for the particular application of interest.

\begin{table}[ht]
\centering
\begin{tabular}{cccr} \hline
Model & \multicolumn{1}{c}{Constraints on $\bc$} & $\bc^*$ & \multicolumn{1}{c}{$\rho_a$} \\ \cmidrule(l){1-1}\cmidrule(l){2-2}\cmidrule(l){3-3}\cmidrule(l){4-4}
1 & $2c_2 = 21c_4$, $c_1 \le 3c_4$, $c_3 \ge 12c_4$ & $[0.077, 0.269, 0.308, 0.026, 0.320]$ & $\mathbf{0.426}$ \\
2 & $2c_2 = 21c_4$, $c_1 \le 3c_4$ & $[0.171, 0.6, 0.001, 0.057, 0.171]$ & $\mathbf{0.846}$ \\
3 & $4c_2 = 21c_4$, $c_1 \le 3c_4$  & $[0.245, 0.429, 0.031, 0.082, 0.214]$ & $\mathbf{0.906}$\\
4 & None & $[0.0002, 0.0001, 0.0001, 0.9995, 0.0001]$ & $\mathbf{0.999}$ \\
\hline
\end{tabular}
\caption{Cost estimates from different inverse APP models without scaling.}
\label{tab:APP_results_nonscaled}
\end{table}

\hfill
\subsection{Forward and Inverse Formulations for Cancer Therapy}\label{appendix:IMRT}
In this section we provide the full formulations for the forward and inverse IMRT treatment plan optimization problem.  We refer the interested reader to~\citet{Chan14} for additional details. Let $\mathcal{B}$ be the set of beamlets and $w_b$ be the intensity delivered by beamlet $b\in\mathcal{B}$.  A patient's anatomy is discretized into volume elements called voxels.  We denote by $D_{v,b}$ the dose deposited to voxel $v$ from unit intensity of beamlet $b$.  Let $\mathcal{K}$ be the set of all objectives.  For any $k\in\mathcal{K}$, let $\mathcal{O}_k$ be the set of voxels in the OAR associated with objective $k$.  We also let $\mathcal{V}$ be the set of voxels in the anatomy.  Lastly, let $\alpha_k \ge 0$ denote the weight assigned to objective $k$.   The forward multi-objective IMRT treatment plan optimization problem is given by
\begin{equation}\label{eq:IMRT_full}
\begin{alignedat}{1}
\underset{\bw}{\text{minimize}}       \quad & \sum_{k\in\mathcal{K}} \alpha_k f_k(\bw)\\
\text{subject to}  \quad & \sum_{b\in\mathcal{B}}D_{v,b}w_b \ge \ell_v,  \quad \forall v \in \mathcal{T}, \\
                         & \sum_{b\in\mathcal{B}}D_{v,b}w_b \le u_v,  \quad \forall v \in \mathcal{V},\\
                         & \frac{\beta_1}{|\mathcal{B}|}\underset{b'\in\mathcal{B}}{\sum}w_{b'} \leq\ w_b \le \frac{\beta_2}{|\mathcal{B}|}\underset{b'\in\mathcal{B}}{\sum}w_{b'}, \quad \forall b \in \mathcal{B},\\
                         & w_b \ge 0, \quad \forall b \in \mathcal{B},
\end{alignedat}
\end{equation}
where organ $i \in \mathcal{I} \subseteq \mathcal{K}$ is associated with objective function type (a):
\begin{equation}\label{eq:obj(a)}
f_i(\bw) =  \frac{1}{|\mathcal{O}_i|}\sum_{v\in \mathcal{O}_{i}} \max\left\{0, \sum_{b\in\mathcal{B}}D_{v,b}w_b - \theta^i_v \right\},
\end{equation}
and organ $j \in \mathcal{J} \subseteq \mathcal{K}$ is associated with objective function type (b):
\begin{equation}\label{eq:obj(b)}
f_j(\bw) = \max_{v \in \mathcal{O}_{j}} \left\{\sum_{b\in \mathcal{B}} D_{v,b}w_b\right\}.
\end{equation}

By introducing auxiliary variables $\xi$ and $\bz$ corresponding to the above-threshold doses (objective type (a)) and max dose (objective type (b)), respectively, formulation~\eqref{eq:IMRT_full} can be written as an LP:
\begin{equation}\label{eq:IMRT_full2}
\begin{alignedat}{1}
\underset{\bw,\bxi,\bz}{\text{minimize}}       \quad & \sum_{i\in\mathcal{I}}\frac{\alpha_i}{|\mathcal{O}_i|} \sum_{v\in \mathcal{O}_i} \xi_{i,v} + \sum_{j\in\mathcal{J}} \alpha_j z_j\\
\text{subject to}  \quad &\xi_{i,v} \ge \sum_{b\in\mathcal{B}}D_{v,b}w_b - \theta_v^i, \quad \forall v \in \mathcal{O}_i, i \in \mathcal{I}, \\
                         & z_j \ge \sum_{b\in\mathcal{B}}D_{v,b}w_b, \quad \forall v \in \mathcal{O}_j, j \in \mathcal{J}, \\
                         & \sum_{b\in\mathcal{B}}D_{v,b}w_b \ge \ell_v,  \quad \forall v \in \mathcal{T}, \\
                         & \sum_{b\in\mathcal{B}}D_{v,b}w_b \le u_v,  \quad \forall v \in \mathcal{V}, \\
                         & \frac{\beta_1}{|\mathcal{B}|}\underset{b'\in\mathcal{B}}{\sum}w_{b'} \leq\ w_b \le \frac{\beta_2}{|\mathcal{B}|}\underset{b'\in\mathcal{B}}{\sum}w_{b'}, \quad \forall b \in \mathcal{B},\\
                         & \xi_{i,v} \ge 0, \quad \forall v \in \mathcal{O}_i, i \in \mathcal{I},\\
                         & z_j \ge 0, \quad j \in \mathcal{J},\\
                         & w_b \ge 0, \quad \forall b \in \mathcal{B}.
\end{alignedat}
\end{equation}
Let $\bp, \bq, \br, \bs, \bt^1$ and $\bt^2$ be the dual variables associated with the constraints in formulation~\eqref{eq:IMRT_full2} ($\bt^1$ and $\bt^2$ for the left and right inequalities in the fifth constraint, respectively).   Given a dose distribution $\bD\hat\bw$ from a clinical treatment plan, the corresponding $\hat\bxi$ and $\hat\bz$  values can be calculated. The relative duality gap inverse formulation is:

\begin{equation}\label{eq:IMRT_full_inverse}
\begin{alignedat}{1}
\underset{\bp,\,\bq,\,\br,\,\bs,\,\bt^1,\,\bt^2,\,\balpha}{\textrm{minimize}}\quad & \sum_{i\in\mathcal{I}}\frac{\alpha_i}{|\mathcal{O}_i|} \sum_{v\in \mathcal{O}_i} \hat\xi_{v,i} + \sum_{j\in\mathcal{J}} \alpha_j \hat{z}_j\\
    \textrm{subject to}\quad & -\sum_{i\in\mathcal{I}}\sum_{v\in \mathcal{O}_i}D_{v,b}p_{v,i}-\sum_{j\in\mathcal{J}}\sum_{v\in\mathcal{O}_j}D_{v,b}q_{v,j}+\sum_{v\in\mathcal{T}}D_{v,b}r_{v}-\sum_{v\in \mathcal{V}}D_{v,b}s_{v}\\
    &\qquad+t^1_b- \frac{\beta_1}{|\mathcal{B}|}\sum_{b' \in \mathcal{B}}t^1_{b'} -t^2_b+\frac{\beta_2}{|\mathcal{B}|}\sum_{b' \in \mathcal{B}}t^2_{b'} \leq 0, \quad \forall b \in B,\\
&p_{v,i} \leq \frac{\alpha_i}{|\mathcal{O}_i|},\quad\forall v\in \mathcal{O}_i,\forall i\in \mathcal{I},\\
&\sum_{v\in\mathcal{O}_j}q_{v,j} = \alpha_j,\quad\forall j\in \mathcal{J},\\
&-\sum_{i \in \mathcal{I}}\sum_{v\in \mathcal{O}_i}\theta_v^i p_{v,i}+\sum_{v\in\mathcal{T}}\ell_{v}r_{v}-\sum_{v\in\mathcal{V}}u_v s_{v} = 1, \\
&\alpha_k \geq0,\quad\forall k\in \mathcal{K},\\
&p_{v,i} \geq0,\quad\forall v\in \mathcal{O}_i,\forall i\in\mathcal{I},\\
&q_{v,j} \geq0,\quad\forall v\in \mathcal{O}_j,\forall j\in\mathcal{J},\\
&r_v \geq 0,\quad\forall v\in \mathcal{T}, \quad s_v \geq 0,\quad\forall v\in\mathcal{V},\\
&t^1_b \geq 0,\quad\forall b\in \mathcal{B}, \quad t^2_b \geq 0,\quad\forall b\in \mathcal{B}.
\end{alignedat}
\end{equation}
Note that the structure of \eqref{eq:IMRT_full_inverse} differs slightly from the general relative duality gap inverse model~\eqref{eq:GIO_r_formulation}. The difference is due to the fact that $\balpha$ is assumed to be nonnegative and $f_k(\bw)\ge 0$ for all $k\in\mathcal{K}$ for every feasible $\bw$. In this case, the objective function simplifies to $\epsilon_r$, which can be minimized by fixing the denominator and minimizing the numerator, followed by a post-hoc re-normalization (cf. formulation (4) from~\citet{Chan14}). See Proposition~\ref{GIO_r_LP_alternative} and formulation~\ref{eq:GIO_r_single}.

\subsection{Proofs}\label{sec:proofs}

\proof{Proof of Proposition~\ref{thm:IO_feasibility}}
($\Rightarrow$) We condition on the feasibility of $\hatX$. If $\hatX$ is feasible, then feasibility of $\mathbf{IO}(\hatX)$ implies the optimality conditions are met, so $\hatX$ is optimal. Now suppose $\hatX$ is infeasible and, to derive a contradiction, satisfies $\ba_i'\hatX < b_i$ for all $i \in I$.  Feasibility of $\mathbf{IO}(\hatX)$ forces $\by'(\bA\hatX-\bbb)=\bzero$, which can be satisfied only if $\by = \bzero$, violating the normalization constraint.

($\Leftarrow$) If $\hatX \in \bX^{\textrm{OPT}}$, then dual feasibility and strong duality hold with non-zero $\by$ and $\bc$.  Thus, $\by$ and $\bc$, both normalized by $\|\bc\|_1$, are feasible for $\mathbf{IO}(\hatX)$.  If $\hatX \in \{\bx \not\in \bX \;|\; \ba_i'\bx \ge b_i \textrm{ for some } i \in I\}$, there exists $i \in I$ such that $\ba_i' \hatX \ge b_i$ and $k \in I$ such that $\ba_{k}' \hatX < b_{k}$.  It is easily checked that $\by = \beee_i - (\ba_i'\hatX - b_i)/(\ba_k'\hatX-b_k)\beee_k$ and $\bc = y_i\ba_i + y_k\ba_k$, both normalized by $\|\bc\|_1$, are feasible for $\mathbf{IO}(\hatX)$. $\Box$ \\
\endproof

\proof{Proof of Proposition~\ref{thm:GIO_feasibility}}
Since the objective has a lower bound of zero, it suffices to show that $\mathbf{GIO}(\hatX)$ is feasible to prove the first statement. If $\hatX \in \bX^{\textrm{OPT}}$, clearly $\mathbf{GIO}(\hatX)$ is feasible (optimal) with $\bepsilon = \bzero$.  If $\hatX \in \bX\backslash\bX^{\textrm{OPT}}$, there exists a dual feasible $\by$, appropriately scaled, such that $\bA'\by = \bc, \|\bc\|_1 = 1$, and $\bbb'\by < \bc'\hatX$.  Since $\bepsilon$ is unrestricted, it can be chosen so that $\bc'\bepsilon = \bc'\hatX - \bbb'\by$. Thus, $\mathbf{GIO}(\hatX)$ is feasible and has an optimal solution. The proof of the second statement is trivial, observing the correspondence between the strong duality constraints in $\mathbf{IO}(\hatX)$ and $\mathbf{GIO}(\hatX)$. $\Box$ \\
\endproof

The proof of Theorem~\ref{thm:GGIO} requires showing that an optimal solution to $\mathbf{GIO}(\hatX)$ can be found by 1) projecting $\hatX$ to one of the hyperplanes defining $\bX$, and 2) characterizing the point on to which $\hatX$ is projected. For the latter step, we leverage Theorem \ref{thm:Mangasarian} below, which characterizes projections of a point on to a hyperplane under general norms.

\begin{theorem}\emph{[\citet[Thm. 2.1]{Mangasarian99}]}\label{thm:Mangasarian}
Let $\| \cdot \|$ be a norm defined over $\mathbb{R}^n$
and let $\mathcal{H}= \{\bx \in \mathbb{R}^n\;|\; \ba'\bx = b\}$, $\ba \ne \bzero$, $\ba \in \mathbb{R}^n$, $b \in \mathbb{R}$.
Let $\bx^0 \in \mathbb{R}^n\setminus\mathcal{H}$. A projection $\pi(\bx^0) \in \mathcal{H}$ under $\| \cdot \|$ is given by $\pi(\bx^0) = \bx^0 - {(\ba'\bx^0 - b)\bv(\ba)}/{\|\ba\|^{D}}$, where $\bv(\ba) \in {\arg\max}_{\|\bv\|=1} \ba'\bv$, and $\|\bx^0 - \pi(\bx^0)\| = |\ba'\bx^0 - b|/\|\ba\|^D$.
\end{theorem}

\proof{Proof of Theorem~\ref{thm:GGIO}}
Let $(\yinv, \cinv, \bepsilon^*)$ be an optimal solution to $\mathbf{GIO}(\hatX)$, and define $\bx^* = \hatX - \bepsilon^*$.  By Proposition~\ref{thm:GIO_feasibility}, $(\yinv, \cinv)$ is a feasible solution to $\mathbf{IO}(\bx^*)$, which, by Proposition~\ref{thm:IO_feasibility}, implies that either $\bx^* \in \boundary(\cinv)$ or $\bx^* \in \{\bx \not\in \bX \;|\; \ba_i'\bx \ge b_i \textrm{ for some } i \in I\}$. We will show the former is true by showing that $\bx^*$ is primal feasible. Suppose to the contrary that $\bx^*$ is infeasible. Then there exist a $j$ such that $\ba_j'(\hatX-\bepsilon^*) < b_j$ and a primal feasible point $\hat\bx := \hatX - \hat\bepsilon = \lambda\hatX + (1-\lambda)(\hatX-\bepsilon^*)$, $\lambda \in (0,1)$, such that $\ba_j'\hat\bx = b_j$. Since $\hat\bx$ is on the boundary of $\bX$, i.e., $\hat\bx\in\boundary$, $\mathbf{IO}(\hat\bx)$ is feasible by Proposition~\ref{thm:IO_feasibility}. If we let $(\hat\by, \hat\bc)$ be a feasible solution to $\mathbf{IO}(\hat\bx)$, then by Proposition~\ref{thm:GIO_feasibility}, $(\hat\by,\hat\bc, \hat\bepsilon)$ is feasible for $\mathbf{GIO}(\hatX)$, where $\hat\bepsilon = (1-\lambda)\bepsilon^*$. However, $\|\hat\bepsilon\|_L = \|(1-\lambda)\bepsilon^*\|_L < \|\bepsilon^*\|_L$, which contradicts the optimality of $(\by^*, \bc^*, \bepsilon^*)$. Thus, it must be that $\bx^* \in \boundary(\cinv)$. In particular, $\bx^*$ is feasible and there exists $i$ such that $\ba_i'\bx^* = b_i$.

By optimality of $\bepsilon^*$, $||\bepsilon^*||_L=||\bx^*-\hatX||_L$ is the minimum distance from $\hatX$ to $\bX^{\textrm{OPT}}$, i.e., the minimum of all distances from $\hatX$ to the hyperplanes defining $\bX$. Applying Theorem 2.1 of~\citet{Mangasarian99} for each $i \in I$, $\|\bepsilon^*\|_L = \underset{i\in I} \min \{ {(\ba_i'\hatX - b_i)}/{\|\ba_i\|_L^D} \}$. Let $\istar \in \underset{i \in I}{\arg\min} \{ {(\ba_i'\hatX - b_i)}/{\|\ba_i\|_L^D} \}$. Then $\invoptX$ satisfies $\ba_{\istar}'\invoptX = b_{\istar}$. Consequently, $\bc^*$ can be written as $\lambda \ba_{i^*}$ for some $\lambda > 0$. The constraints $\|\bc\|_1 = 1$ and $\bA'\by = \bc$ imply $\cinv = {\ba_{\istar}}/{\|\ba_{\istar}\|_1}$ and $\by^* = \beee_{\istar}/{\|\ba_{\istar}\|_1}$. The expression for $\bepsilon^*$ follows directly from Theorem 2.1 of~\citet{Mangasarian99}. $\Box$ \\
\endproof

\proof{Proof of Proposition~\ref{thm:equivalenceGIOa}}
First, the objective function of $\mathbf{GIO}(\hatX)$ with the chosen infinity norm and structure of $\bepsilon$ matches the objective function of $\mathbf{GIO}_a(\hatX)$: $\| \bepsilon \|_{\infty} = \| \epsilon_a \sgn(\bc) \|_{\infty} = |\epsilon_a|\; \|\sgn(\bc)\|_{\infty} = \epsilon_a.$ Finally, the strong duality constraint of $\mathbf{GIO}(\hatX)$ becomes: $\bbb'\by = \bc'(\hatX-\bepsilon) = \bc'\left(\hatX - \epsilon_a \sgn(\bc)\right) =  \bc'\hatX - \epsilon_a$, since $ \bc'\sgn(\bc) = \|\bc\|_1 = 1$. $\Box$ \\
\endproof

\proof{Proof of Proposition~\ref{thm:GIOa_reformulation}}
By Proposition~\ref{thm:equivalenceGIOa}, $\mathbf{GIO}_{\infty}(\hatX)$ (i.e., $\mathbf{GIO}(\hatX)$ with $\|\cdot\|_L = \|\cdot\|_{\infty}$) is a relaxed version $\mathbf{GIO}_a(\hatX)$ with the constraint $\bepsilon = \epsilon_a \sgn(\bc)$ omitted.  Applying Theorem 1, $(\by^*, \bc^*, \bepsilon^*) = (\beee_{i^*}/\|\ba_{i^*}\|_1, \ba_{i^*}/\|\ba_{i^*}\|_1, (\ba_{\istar}'\hatX - b_{\istar}) \sgn(\ba_{i^*}) /\|\ba_{\istar}\|_1)$ is an optimal solution to $\mathbf{GIO}_{\infty}(\hatX)$, where $i^* \in \arg\min_{i \in I} (\ba_{i}'\hatX - b_{i}) /\|\ba_{i}\|_1$. This solution is feasible for $\mathbf{GIO}_a(\hatX)$ since there always exists an $\epsilon_a$ such that $\bepsilon^* = \epsilon_a \sgn(\bc^*)$. In particular, $\epsilon_a = (\ba_{\istar}'\hatX - b_{\istar})/\|\ba_{\istar}\|_1$. Since $(\by^*, \bc^*, \bepsilon^*)$ is optimal for $\mathbf{GIO}_{\infty}(\hatX)$ and feasible for $\mathbf{GIO}_a(\hatX)$, it must also be optimal for $\mathbf{GIO}_a(\hatX)$, and hence $\epsilon^*_a = (\ba_{\istar}'\hatX - b_{\istar})/\|\ba_{\istar}\|_1$. $\Box$ \\
\endproof

\proof{Proof of Proposition~\ref{thm:equivalenceGIOr}}
First, the objective function of $\mathbf{GIO}(\hatX)$ with the chosen weighted infinity norm and structure of $\bepsilon$ matches the objective function of $\mathbf{GIO}_r(\hatX)$: $\| \bepsilon \|_{\infty, 1/|\bbb'\by|} = \|  \bbb'\by(\epsilon_r - 1) \sgn(\bc) \|_{\infty}/|\bbb'\by| = |\bbb'\by(\epsilon_r - 1)| \|\sgn(\bc)\|_{\infty}/|\bbb'\by| = |\epsilon_r - 1|$. Second, the strong duality constraint of $\mathbf{GIO}(\hatX)$ becomes:
$\bc'\hatX = \bbb'\by + \bc'\bepsilon = \bbb'\by + \bc'(\bbb'\by(\epsilon_r - 1)\sgn(\bc)) = \bbb'\by + \bbb'\by(\epsilon_r - 1)  = \epsilon_r \bbb'\by$, where the second last equality is due to $\bc'\sgn(\bc) = 1$. $\Box$ \\
\endproof

\proof{Proof of Proposition~\ref{thm:GIOr_reformulation}}
By Proposition~\ref{thm:equivalenceGIOr}, $\mathbf{GIO}_{\infty,K}(\hatX)$ (i.e., $\mathbf{GIO}(\hatX)$ with $\|\cdot\|_L = \|\cdot\|_{\infty,K}$) where $K = 1/|\bbb'\by|$ is a relaxed version of $\mathbf{GIO}_r(\hatX)$ with the constraint $\bepsilon = \bbb'\by(\epsilon_r - 1) \sgn(\bc)$ omitted. To apply Theorem~\ref{thm:GGIO} to $\mathbf{GIO}_{\infty,K}(\hatX)$, we first need to derive the dual of the weighted infinity norm $\|\cdot\|_{\infty,K}$. We have $\| \ba_i \|^D_{\infty, 1/|\bbb'\by|} = \underset{\bx} \sup \{ \ba_i'\bx : \|\bx\|_{\infty, 1/|\bbb'\by|} \le 1 \} = \underset{\bx} \sup \{ \ba_i'\bx : \|\bx\|_{\infty} /|\bbb'\by| \le 1 \} = \underset{\bx} \sup \{ \ba_i'\bx : \|\bx\|_{\infty} \le |\bbb'\by| \} \le \underset{\bx} \sup \{\|\ba_i\|_1 \|\bx\|_{\infty} : \|\bx\|_{\infty} \le |\bbb'\by|\} \le \|\ba_i\|_1 |\bbb'\by|,$ where the first inequality is due to H{\"o}lder's inequality. If we pick $\bx = |\bbb'\by|\sgn(\ba_i)$, we achieve equality throughout. Hence, $\| \ba_i \|^D_{\infty, 1/|\bbb'\by|}  = \|\ba_i\|_1 |\bbb'\by|$.

By definition, the choice of $\bx$ above that achieves equality defines $\bv(\ba_i)$. So, $\bv(\ba_i) = |\bbb'\by|\sgn(\ba_i)$. Applying Theorem~\ref{thm:GGIO}, $(\by^*, \bc^*, \bepsilon^*)
= \left(\beee_{i}/\|\ba_{i}\|_1, \ba_{i}/\|\ba_{i}\|_1, (\ba_{i}'\hatX - b_{i})\bv(\ba_{i})/\|\ba_{i}\|_{L}^{D} \right)
= \left(\beee_{i}/\|\ba_{i}\|_1, \ba_{i}/\|\ba_{i}\|_1, (\ba_{i}'\hatX - b_{i})\sgn(\ba_{i})/\|\ba_{i}\|_1\right)$ for some $i$. Given the optimal solution structure, the optimal objective value is $\|\bepsilon^*\|_{\infty} / |\bbb'\by^*| = \|(\ba_{i^*}'\bx^0-b_{i^*})\sgn(\ba_{i^*})/\|\ba_{i^*}\|_1\|_{\infty} / (|b_{i^*}|/\|\ba_{i^*}\|_1) = |(\ba_{i^*}'\bx^0 - b_{i^*})/b_{i^*}| = (\ba_{i^*}'\bx^0 - b_{i^*})/|b_{i^*}|$ where $i^* \in \argmin_{i} \{(\ba_{i}' \bx^0 - b_i)/|b_i | \}$.

Finally, $\epsilon_r^* = \bepsilon^*/(\sgn(\bc^*)\bbb'\by^*) + 1 = \left((\ba_{i^*}'\bx^0-b_{i^*})\sgn(\ba_{i^*})/\|\ba_{i^*}\|_1\right) / (\sgn(\ba_{i^*})b_{i^*}/\|\ba_{i^*}\|_1) + 1 = (\ba_{i^*}'\bx^0-b_{i^*})/b_{i^*} + 1 = \ba_{i^*}'\bx^0/b_{i^*}$, as required. $\Box$ \\
\endproof

\proof{Proof of Proposition~\ref{GIO_r_LP_alternative}}
There is a one-to-one correspondence between feasible solutions to~\eqref{eq:GIO_r_formulation} and~\eqref{eq:GIO_r_alternative} since they only differ by a scaling constraint. The term $\epsilon_r$ in the objective equals $\bc'\hatX / \bbb'\by$ and is invariant to simultaneous scaling of the $\bc$ and $\by$ vectors. Thus, $(\hat{\by}, \hat{\bc}, \hat{\bepsilon})$ to optimal to~\eqref{eq:GIO_r_alternative} if and only if $(\hat{\by}/\|\hat{\bc}\|_1,\hat{\bc}/\|\hat{\bc}\|_1, \hat{\epsilon}_r)$ is optimal to~\eqref{eq:GIO_r_formulation}. $\Box$ \\
\endproof

\proof{Proof of Theorem~\ref{thm:rho_properties}}
\hfill
\begin{enumerate}
\item Since an optimal $\bepsilon^*$ to $\mathbf{GIO}(\hatX)$ minimizes $\|\bepsilon\|_L$ and the sum $\sum_{i=1}^m \|\bepsilon^i\|_L$ depends only on $\hatX$ and the primal feasible region, $\rho$ is maximized by $\bepsilon^*$. 

\item By definition, an optimal $\bepsilon^*$ to $\mathbf{GIO}(\hatX)$ satisfies $0 \le \|\bepsilon^*\|_L \le \|\bepsilon^i\|_L$ for all $i = 1, \ldots, m$.  Thus, $0 \le \|\bepsilon^*\|_L \le \frac{1}{m}\sum_{i=1}^m \|\bepsilon^i\|_L$, which implies $0 \le \rho \le 1$ if $\frac{1}{m}\sum_{i=1}^m \|\bepsilon^i\|_L \neq 0$. If $\frac{1}{m}\sum_{i=1}^m \|\bepsilon^i\|_L = 0$, we define $\rho := 1$.

\item Let $(\by^{(k)*}, \bc^{(k)*}, \bepsilon^{(k)*})$ be an optimal solution to $\mathbf{GIO}^{(k)}(\hatX)$. This solution is also feasible (but not necessarily optimal) for $\mathbf{GIO}^{(k+1)}(\hatX)$. Thus, $\|\bepsilon^{(k)*}\|_L \ge \|\bepsilon^{(k+1)*}\|_L$ and the result follows. 

\item Let $\bx_1$ and $\bx_2$ be distinct feasible solutions that, without loss of generality, satisfy $\rho(\bx_1) \le \rho(\bx_2)$.  We wish to show that $\rho(\bar\bx) \le \max\left\{\rho(\bx_1),\rho(\bx_2)\right\} = \rho(\bx_2)$, where $\bar\bx=\lambda\bx_1+(1-\lambda)\bx_2, \lambda \in (0,1)$, or equivalently $1 - \rho(\bar\bx) \ge \min \{1 - \rho(\bx_1), 1 - \rho(\bx_2)\} = 1 - \rho(\bx_2).$  Let
\begin{equation}\label{eq:rho_bar}
1 - \rho(\bar\bx) = \frac{\|\bar\bepsilon^*\|_L}{\frac{1}{m}\sum_{i=1}^m \|\bar\bepsilon^{i}\|_L},
\end{equation}
and
\begin{equation}
1 - \rho(\bx_k) = \frac{\|\bepsilon_k^*\|_L}{\frac{1}{m}\sum_{i=1}^m \|\bepsilon_k^i\|_L},
\end{equation}
for $k = 1, 2$, where $\bar\bepsilon^i$ and $\bepsilon_k^i$ are optimal solutions to problem~\eqref{eq:eps_p^i} with respect to $\bar\bx$ and $\bx_k$, $k = 1, 2$, respectively.  Let $\bar i^* \in \underset{i \in I} {\text{argmin}} \{ (\ba_i'\bar\bx - b_i)/\|\ba_i\|_L^D\}$, $i_1^* \in \underset{i \in I} {\text{argmin}} \{ (\ba_i'\bx_1 - b_i)/\|\ba_i\|_L^D\}$, and $i_2^* \in \underset{i \in I} {\text{argmin}} \{ (\ba_i'\bx_2 - b_i)/\|\ba_i\|_L^D\}$. First, we bound the numerator of~\eqref{eq:rho_bar}:
\begin{equation}\label{eq:rho_num}
\begin{alignedat}{1}
\|\bar\bepsilon^*\|_L & = \frac{{\ba_{\bar{i}^{*}}}'\bar\bx - b_{\bar{i}^{*}}}{\|\ba_{\bar{i}^{*}}\|_L^D} = \frac{{\ba_{\bar{i}^{*}}}'(\lambda\bx_1 + (1-\lambda)\bx_2) - b_{\bar{i}^{*}}}{\|\ba_{\bar{i}^{*}}\|_L^D} = \lambda\frac{{\ba_{\bar{i}^{*}}}'\bx_1- b_{\bar{i}^{*}}}{\|\ba_{\bar{i}^{*}}\|_L^D} + (1-\lambda)\frac {{\ba_{\bar{i}^{*}}}'\bx_2- b_{\bar{i}^{*}}}{\|\ba_{\bar{i}^{*}}\|_L^D} \\
& \ge \lambda\frac{{\ba_{i_1^*}}'\bx_1- b_{i_1^*}}{\|\ba_{i_1^*}\|_L^D} + (1-\lambda)\frac {{\ba_{i_2^*}}'\bx_2- b_{i_2^*}}{\|\ba_{i_2^*}\|_L^D} = \lambda\|\bepsilon_1^*\|_L + (1-\lambda)\|\bepsilon_2^*\|_L,
\end{alignedat}
\end{equation}
where the first and last equalities are from Theorem~\ref{thm:GGIO}, and the inequality is from optimality of $\mathbf{GIO}(\bx_k)$ for $k = 1, 2$.

Next we bound the denominator of~\eqref{eq:rho_bar}. Note that if $\bepsilon_k$ satisfies $\bA(\bx_k-\bepsilon_k) \ge \bbb$, for $k = 1, 2$, then $\bA(\bar\bx-\bepsilon) \ge \bbb$ if $\bepsilon = \lambda\bepsilon_1 + (1-\lambda)\bepsilon_2$. Similarly, if $\ba_i'(\bx_k-\bepsilon_k) = b_i$ for $k=1,2$, then $\ba_i'(\bar\bx-\bepsilon) = b_i$. In other words, given feasible solutions $\bepsilon_k$ to problem~\eqref{eq:eps_p^i} with respect to $\bx_k$, $k=1,2$, we can construct a feasible solution $\bepsilon$ to problem~\eqref{eq:eps_p^i} with respect to $\bar\bx$. Thus,
\begin{equation}
\|\bar\bepsilon^i\|_L \le \|\lambda\bepsilon^i_1 + (1-\lambda)\bepsilon^i_2\|_L \le \lambda\|\bepsilon^i_1\|_L + (1-\lambda)\|\bepsilon^i_2\|_L,
\end{equation}
which implies
\begin{equation}~\label{eq:rho_demon}
\frac{1}{m}\sum_{i=1}^m\|\bar\bepsilon^i\|_L \le \lambda\left(\frac{1}{m}\sum_{i=1}^m\|\bepsilon_1^i\|_L\right) + (1 - \lambda)\left(\frac{1}{m}\sum_{i=1}^m\|\bepsilon_2^i\|_L\right).
\end{equation}
Putting~\eqref{eq:rho_num} and~\eqref{eq:rho_demon} together,
\begin{equation}~\label{eq:rho_combined}
\frac{\|\bar\bepsilon^{*}\|_L}{\frac{1}{m}\sum_{i=1}^m \|\bar\bepsilon^{i}\|_L} \ge   \frac{\lambda\|\bepsilon_1^*\|_L + (1-\lambda)\|\bepsilon_2^*\|_L}{\lambda(\frac{1}{m}\sum_{i=1}^m\|\bepsilon_1^i\|_L) + (1 - \lambda)(\frac{1}{m}\sum_{i=1}^m\|\bepsilon_2^i\|_L)}.
\end{equation}
What remains is to show that
\begin{equation}\label{eq:rho_remains}
\frac{\lambda\|\bepsilon_1^*\|_L + (1-\lambda)\|\bepsilon_2^*\|_L}{\lambda(\frac{1}{m}\sum_{i=1}^m\|\bepsilon_1^i\|_L) + (1 - \lambda)(\frac{1}{m}\sum_{i=1}^m\|\bepsilon_2^i\|_L)} \ge \frac{\|\bepsilon^{*}_2\|_L}{\frac{1}{m}\sum_{i=1}^m \|\bepsilon^{i}_2\|_L},
\end{equation}
which follows from straightforward algebraic manipulation, along with the original assumption that
\begin{equation}
1-\rho(\bx_1) = \frac{\|\bepsilon^{*}_1\|_L}{\frac{1}{m}\sum_{i=1}^m \|\bepsilon^{i}_1\|_L} \ge \frac{\|\bepsilon^{*}_2\|_L}{\frac{1}{m}\sum_{i=1}^m \|\bepsilon^{i}_2\|_L} = 1-\rho(\bx_2).
\end{equation}
Thus, combining~\eqref{eq:rho_combined} and~\eqref{eq:rho_remains}, we have $1-\rho(\bar\bx) \ge 1-\rho(\bx_2)$, as desired. $\Box$ \\
\end{enumerate}
\endproof

\proof{Proof of Proposition~\ref{thm:distance_quasiconvexity}}
Suppose $\|\bx_1 - \bx^*\|_L \le \|\bx_2 - \bx^*\|_L$. If $\|\bx_1 - \bx^*\|_L = \|\bx_2 - \bx^*\|_L$, then $\bx_1 = \bx_2$ and $\rho(\bx_1) = \rho(\bx_2)$.
If $\|\bx_1 - \bx^*\|_L < \|\bx_2 - \bx^*\|_L$ then, using the definition of $\bx_1$ and $\bx_2$, we find that $\|\lambda_1\bx^{\min} + (1 - \lambda_1)\bx^* - \bx^*\|_L < \|\lambda_2\bx^{\min} + (1 - \lambda_2)\bx^*- \bx^*\|_L$, which implies $\lambda_1 < \lambda_2$. Additionally, from the definition of $\bx_1$ it follows that $\bx^* = (\bx_1 - \lambda_1\bx^{\min})/(1 - \lambda_1)$; substituting this expression into the equation for $\bx_2$, we get $\bx_2 = (\lambda_2 - \lambda_1)\bx^{\min}/(1- \lambda_1)$ + $(1 - \lambda_2)\bx_1/(1- \lambda_1)$. Since $\lambda_1 < \lambda_2$, we now know there exists $\lambda \in (0, 1)$ such that $\bx_2 = \lambda \bx^{\min} + (1 - \lambda) \bx_1$. By quasiconvexity of $\rho$, we know that $\rho(\bx_2) \le \max\{\rho(\bx^{\min}),\rho(\bx_1)\} = \rho(\bx_1)$.$\Box$ \\
\endproof

\proof{Proof of Proposition~\ref{thm:approxrho_properties}}
The first three properties are immediate since the only difference between $\tilde\rho$ and $\rho$ is in the denominator, which is independent of the inverse process. For the quasiconvexity proof, let $\bx_1$ and $\bx_2$ be distinct feasible solutions that, without loss of generality, satisfy $\tilde\rho(\bx_1) \le \tilde\rho(\bx_2)$. We wish to show that $\tilde\rho(\bar\bx) \le \max\{\tilde\rho(\bx_1),\tilde\rho(\bx_2)\} = \tilde\rho(\bx_2)$, where $\bar\bx=\lambda\bx_1+(1-\lambda)\bx_2, \lambda \in (0,1)$. We have
\begin{align*}
1-\tilde\rho(\bar\bx) & = \frac{\|\bar\bepsilon^*\|_L}{\frac{1}{m}\sum_{i=1}^m\|\tilde{\bar\bepsilon}^i\|_L}\\
& \ge \frac{\lambda\|\bepsilon_1^*\|_L+(1-\lambda)\|\bepsilon_2^*\|_L}{\frac{1}{m}\sum_{i=1}^m\|\tilde{\bar\bepsilon}^i\|_L}\\
& = \frac{\lambda\|\bepsilon_1^*\|_L+(1-\lambda)\|\bepsilon_2^*\|_L}{\frac{1}{m}\sum_{i=1}^m(\lambda\|\tilde\bepsilon_1^i\|_L+(1-\lambda)\|\tilde\bepsilon_2^i\|_L)}\\
& \ge \frac{\|\bepsilon_2^*\|_L}{\frac{1}{m}\sum_{i=1}^m\|\tilde\bepsilon_2^i\|_L}\\
& = 1-\tilde\rho(\bx_2),
\end{align*}
where the first inequality is due to~\eqref{eq:rho_num}, the second equality is due to the linearity of $\|\tilde{\bar\bepsilon}^i\|_L$ with respect to $\bar\bx$, and the last inequality follows from straightforward algebraic manipulation along with the original assumption that $\tilde\rho(\bx_1) \le \tilde\rho(\bx_2)$. $\Box$ \\
\endproof

\bibliography{inverseOptimizationBibliography}
\bibliographystyle{ormsv080}

\end{document}